\documentclass[aos,MSNbibl,nameyear,seceqn,dvips]{arximspdf}

\doi{10.1214/13-AOS1151} 
\volume{41}
\issue{6}
\pubyear{2013}
\firstpage{2905}
\lastpage{2947}

\makeatletter

\newtheorem{theorem}{Theorem}[section]
\newtheorem{lemma}{Lemma}[section]
\newtheorem{corollary}{Corollary}[section]
\newtheorem{proposition}{Proposition}[section]

\newproclaim{condition}{Condition}
\newproclaim{definition}{Definition}
\newproclaim{remark}{Remark}
\newproclaim{Conjecture}{Conjecture}

\newcommand{\rrVert}{\Vert}
\newcommand{\rrvert}{\vert}
\newcommand{\llVert}{\Vert}
\newcommand{\llvert}{\vert}

\renewcommand{\Pr}{\mathbb{P}}

\def\tr{\operatorname{trace}}

\newcommand{\prob}{\mathbb{P}}
\newcommand{\expect}{\mathbb{E}}
\def\E{\expect}
\def\Var{\operatorname{Var}}

\newcommand{\diag}{\operatorname{diag}}

\def\Real{\mathbb{R}}

\def\normal{\mathcal{N}}

\def\Sphere{\mathbb{S}}

\newcommand{\colspan}{\operatorname{col}}
\newcommand{\vecspan}{\operatorname{span}}

\def\SampleCov{S_n}

\makeatother

\begin{document}
\begin{frontmatter}

\title{Minimax sparse principal subspace estimation in high~dimensions}
\runtitle{Minimax sparse principal subspace estimation}

\begin{aug}
\author[A]{\fnms{Vincent Q.} \snm{Vu}\corref{}\thanksref{t1}\ead[label=e1]{vqv@stat.osu.edu}}
\and
\author[B]{\fnms{Jing} \snm{Lei}\thanksref{t2}\ead[label=e2]{jinglei@andrew.cmu.edu}}
\runauthor{V. Q. Vu and J. Lei}
\affiliation{Ohio State University and Carnegie Mellon University}
\address[A]{Department of Statistics\\
Ohio State University\\
Columbus, Ohio 43210\\
USA\\
\printead{e1}}
\address[B]{Department of Statistics\\
Carnegie Mellon University\\
Pittsburgh, Pennsylvania 15213\\
USA\\
\printead{e2}}
\end{aug}

\thankstext{t1}{Supported in part by NSF Postdoctoral Fellowship DMS-09-03120.}
\thankstext{t2}{Supported in part by NSF Grant BCS-0941518.}

\received{\smonth{2} \syear{2013}}
\revised{\smonth{6} \syear{2013}}

%
\begin{abstract}
We study sparse principal components analysis in high dimensions, where
$p$ (the number of variables) can be much larger than $n$ (the number
of observations), and analyze the problem of estimating the subspace
spanned by the principal eigenvectors of the population covariance
matrix. We introduce two complementary notions of $\ell_q$ subspace
sparsity: row sparsity and column sparsity. We prove nonasymptotic
lower and upper bounds on the minimax subspace estimation error for $0
\leq q \leq1$. The bounds are optimal for row sparse subspaces and
nearly optimal for column sparse subspaces, they apply to general
classes of covariance matrices, and they show that $\ell_q$ constrained
estimates can achieve optimal minimax rates without restrictive spiked
covariance conditions. Interestingly, the form of the rates matches
known results for sparse regression when the effective noise variance
is defined appropriately. Our proof employs a novel variational
$\sin\Theta$ theorem that may be useful in other regularized spectral
estimation problems.
\end{abstract}

%

\end{frontmatter}

\section{Introduction}
\label{secintroduction}

Principal components analysis (PCA) was introduced in the early 20th century
[\citet{Pearson1901,Hotelling1933}] and is arguably the most well
known and
widely used technique for dimension reduction. It is part of the mainstream
statistical repertoire and is routinely used in numerous and diverse
areas of
application. However, contemporary applications often involve much
higher-dimensional data than envisioned by the early developers of PCA.
In such
high-dimensional situations, where the number of variables $p$ is of
the same
order or much larger than the number of observations $n$, serious difficulties
emerge: standard PCA can produce inconsistent estimates of the principal
directions of variation and lead to unreliable conclusions
[\citet{Johnstone2009,Paul2007,Nadler2008}].\looseness=1

The principal directions of variation correspond to the eigenvectors of the
covariance matrix, and in high-dimensions consistent estimation of the
eigenvectors is generally not possible without additional assumptions about
the covariance matrix or its eigenstructure. Much of the recent development
in PCA has focused on methodology that applies the concept of sparsity
to the
estimation of individual eigenvectors [examples
include \citet
{Jolliffe2003,dAspremont2007,Zou2006,Shen2008,Witten2009,Journee2010}].
Theoretical developments on sparsity and PCA include
consistency [\citet{Johnstone2009,Shen2011}], variable selection properties
[\citet{Amini2009}], rates of convergence and minimaxity
[\citet
{Vu2012}], but
have primarily been limited to results about estimation of the leading
eigenvector. Very recently, \citet{Birnbaum2012} established
minimax lower
bounds for the estimation of individual eigenvectors. However, an open problem
that has remained is whether sparse PCA methods can \emph{optimally} estimate
the subspace spanned by the leading eigenvectors, that is, the \emph{principal
subspace} of variation.

The subspace estimation problem is directly connected to dimension reduction
and is important when there may be more than one principal component of
interest. Indeed, typical applications of PCA use the projection onto the
principal subspace to facilitate exploration and inference of important
features of the data. In that case, the assumption that there are distinct
principal directions of variation is mathematically convenient but unnatural:
it avoids the problem of unidentifiability of eigenvectors by imposing an
artifactual choice of principal axes. Dimension reduction by PCA
should emphasize subspaces rather than eigenvectors.


An important conceptual issue in applying sparsity to
principal subspace estimation is that, unlike the case of sparse vectors,
it is not obvious how to formally
define what is meant by a \emph{sparse principal subspace}.
In this article, we present two complementary notions of sparsity based on
$\ell_q$ (pseudo-) norms: row sparsity and column sparsity.
Roughly, a subspace is row sparse if every one of its orthonormal bases
consists of sparse vectors. In the $q=0$ case, this intuitively means that
a row sparse subspace is generated by a small subset of variables, independent
of the choice of basis. A column sparse subspace, on the other hand, is one
which has some orthonormal basis consisting of sparse vectors. This
means that
the choice of basis is crucial; the existence of a sparse basis is
an implicit assumption behind the frequent use of rotation techniques
by practitioners
to help interpret principal components.

In this paper, we study sparse principal subspace estimation in
high-dimensions. We present nonasymptotic minimax lower and upper
bounds for
estimation of both row sparse and column sparse principal subspaces.
Our upper bounds are constructive and apply to a wide class of distributions
and covariance matrices.
In the row sparse case they are optimal up to constant factors,
while in the column sparse case they are nearly optimal.
As an illustration, one consequence of our results is that the
order of the minimax mean squared estimation error of a row sparse
$d$-dimensional principal subspace (for $d \ll p$) is
\[
R_q \biggl(\frac{\sigma^2}{n}(d+\log p) \biggr)^{1-q/2},\qquad
0\leq q \leq1,
\]
where $\sigma^2$ is the effective noise variance (a function of the
eigenvalues of population covariance matrix) and $R_q$ is a
measure of the sparsity in an $\ell_q$ sense defined in Section~\ref{secprelim}.
Our analysis allows $\sigma$, $R_q$, and $d$ to change with $n$ and $p$.
When $q=0$, the rate has a very intuitive explanation. There are $R_0$ variables
active in generating the principal subspace. For each active variable,
we must estimate the corresponding $d$ coordinates of the basis vectors.
Since we do not know in advance which variables are active, we incur an
additional cost of $\log p$ for variable selection.

To our knowledge, the only other work that has considered sparse principal
subspace estimation is that of \citet{Ma2011}. He proposed a
sparse principal
subspace estimator based on iterative thresholding, and derived its
rate of
convergence under a spiked covariance model (where the covariance
matrix is
assumed to be a rank-$d$ perturbation of the identity) similar to that in
\citet{Birnbaum2012}. He showed that it nearly achieves the
optimal rate
when estimating a single eigenvector, but was not able to track its
dependence on the dimension of the principal subspace.

We obtain the minimax upper bounds by analyzing a sparsity constrained
principal subspace estimator and showing that it attains the optimal
error (up
to constant factors). In comparison to most existing works in the
literature, we show that the upper bounds hold without assuming a spiked
covariance model. This spiked covariance assumption seems to be
necessary for
two reasons. The first is that it simplifies analyses and enables the
exploitation of special properties of the multivariate Gaussian distribution.
The second is that it excludes the possibility of the variables having equal
variances. Estimators proposed by \citet{Paul2007}, \citet
{Johnstone2009}, and \citet{Ma2011} require an initial estimate based
on \emph{diagonal thresholding}---screening out
variables with small sample variances. Such an initial estimate will
not work
when the variables have equal variances or have been standardized. The spiked
covariance model excludes that case and, in particular, does not allow
PCA on
correlation matrices.

A key technical ingredient in our analysis of the subspace estimator is a
novel variational form of the Davis--Kahan $\sin\Theta$ theorem
(see Corollary~\ref{lemcurvature-centered}) that may be useful in
other regularized spectral estimation problems. It allows us to bound the
estimation error using some recent advanced results in empirical
process theory, without
Gaussian or spiked covariance assumptions. The minimax lower bounds
follow the
standard Fano method framework [e.g., \citet{Yu1997}], but their
proofs involve nontrivial
constructions of packing sets in the Stiefel manifold.
We develop a generic
technique that allows us to convert global packing sets without orthogonality
constraints into local packing sets in the Stiefel manifold, followed
by a careful
combinatorial analysis on the cardinality of the resulting matrix class.

The remainder of the paper is organized as follows. In the next
section, we
introduce the sparse principal subspace estimation problem and formally describe
our minimax framework and estimator. In Section~\ref{secresults}, we
present our
main conditions and results, and provide a brief discussion about their
consequences and intuition. Section~\ref{secsketchproof} outlines the
key ideas and
main steps of the proof. Section~\ref{secdiscussion} concludes the paper
with discussion of related problems and practical concerns.
Appendices~\ref{secprooflowerbound},~\ref{secproofupperbound}
contain the details in proving the lower and upper bounds. The major
steps in the proofs require some auxiliary lemmas whose proofs we defer to
Appendices~\ref{secauxproofs},~\ref{secempiricalprocessproofs}.

\section{Subspace estimation}
\label{secprelim}

Let $X_1,\ldots, X_n \in\Real^p$ be independent, identically distributed
random vectors with mean $\mu$ and covariance matrix $\Sigma$. To
reduce the
dimension of the $X_i$'s from $p$ down to $d$, PCA looks for $d$ mutually
uncorrelated, linear combinations of the $p$ coordinates of $X_i$ that have
maximal variance. Geometrically, this is equivalent to finding a
$d$-dimensional linear subspace that is closest to the centered random vector
$X_i - \mu$ in a mean squared sense, and it corresponds to the optimization
problem
%
%
\begin{eqnarray}
\label{eqsubspace-pca} &&\mbox{minimize}\quad \E\bigl\llVert(I_{p} -
\Pi_{\mathcal{G}}) (X_{i} - \mu) \bigr\rrVert_2^2
\nonumber\\[-8pt]\\[-8pt]
&&\mbox{subject to}\quad \mathcal{G} \in\mathbb{G}_{p,d},
\nonumber
\end{eqnarray}
where $\mathbb{G}_{p,d}$ is the Grassmann manifold of $d$-dimensional
subspaces of $\Real^p$, $\Pi_{\mathcal{G}}$ is the orthogonal
projector of $\mathcal{G}$,
and $I_{p}$ is the $p \times p$ identity matrix. [For background on
Grassmann and Stiefel manifolds,
see \citet{Edelman1998,Chikuse2003}.]
There is always at least one $d \leq p$ for which (\ref
{eqsubspace-pca}) has a
unique solution. That solution can be determined by the spectral decomposition
%
%
\begin{equation}
\label{eqpop-spectraldecomposition} \Sigma= \sum_{j=1}^p
\lambda_j v_j v_j^T,
\end{equation}
where $\lambda_1 \geq\lambda_2 \geq\cdots\geq\lambda_p \geq0$
are the
eigenvalues of $\Sigma$ and $v_1,\ldots, v_p \in\Real^p$,
orthonormal, are
the associated eigenvectors. If\vadjust{\goodbreak} $\lambda_d > \lambda_{d+1}$, then the
$d$-dimensional \emph{principal subspace} of $\Sigma$ is
%
%
\begin{equation}
\label{eqprincipal-subspace} \mathcal{S}= \vecspan\{v_1,\ldots,v_d\},
\end{equation}
and the orthogonal projector of $\mathcal{S}$ is given by
$\Pi_{\mathcal{S}} = VV^T$,
where $V$ is the $p \times d$ matrix with columns $v_1,\ldots,v_d$.

In practice, $\Sigma$ is unknown, so $\mathcal{S}$ must be
estimated from
the data. Standard PCA replaces (\ref{eqsubspace-pca}) with an
empirical version.
This leads to the spectral decomposition of the sample covariance
matrix
\[
\SampleCov= \frac{1}{n} \sum_{i=1}^n
(X_i - \bar{X}) (X_i - \bar{X})^T,
\]
where $\bar{X}$ is the sample mean, and estimating
$\mathcal{S}$ by the span
of the leading $d$ eigenvectors of $\SampleCov$. In high-dimensions however,
the eigenvectors of $\SampleCov$ can be inconsistent estimators of the
eigenvectors of $\Sigma$. Additional structural constraints are
necessary for
consistent estimation of $\mathcal{S}$.

\subsection{Subspace sparsity}
The notion of sparsity is appealing and has been used successfully in the
context of estimating vector valued parameters such as the leading eigenvector
in PCA. Extending this notion to subspaces requires care because
sparsity is
inherently a coordinate-dependent concept while subspaces are
coordinate-independent. For a given $d$-dimensional subspace
$\mathcal{G} \in\mathbb{G}_{p,d}$, the set of orthonormal matrices whose
columns span $\mathcal{G}$ is a subset of the Stiefel manifold
$\mathbb{V}_{p,d}$ of $p \times d$ orthonormal matrices.
We will consider two complementary notions of subspace sparsity defined
in terms of those
orthonormal matrices: \emph{row sparsity} and \emph{column sparsity}.

Define the $(2,q)$-norm, $q \in[0,\infty]$, of a $p\times d$ matrix
$A$ as
the usual $\ell_q$ norm of the vector of row-wise $\ell_2$ norms of $A$:
\[
\llVert A \rrVert_{2,q}:= \bigl\llVert\bigl(\matrix{ \llVert
a_{1*} \rrVert_2 & \cdots& \llVert a_{p*}
\rrVert_2 }\bigr) \bigr\rrVert_q,
\]
where $a_{j*}$ denotes the $j$th row of $A$.
(To be precise, this is actually a pseudonorm when $q < 1$.)
Note that $\llVert\cdot\rrVert_{2,q}$
is coordinate-independent, because $\llVert A O \rrVert
_{2,q} = \llVert A \rrVert_{2,q}$ for any
orthogonal matrix $O \in\Real^{d \times d}$. We define the \emph{row sparse
subspaces} using this norm. Let $\colspan(U)$ denotes the span of the
columns of
$U$.
\begin{definition*}[(Row sparse subspaces)]
For $0 \leq q < 2$ and $d \leq R_q \leq d^{q/2} \*p^{1-q/2}$,
\[
\mathcal{M}_{q} (R_q):= \cases{ \bigl\{ \colspan(U)\dvtx U
\in\mathbb{V}_{p,d} \mbox{ and } \llVert U \rrVert_{2,q}^q
\leq R_q \bigr\}, &\quad if $0 < q < 2$\quad and
\vspace*{2pt}\cr
\bigl\{ \colspan(U)\dvtx U
\in\mathbb{V}_{p,d} \mbox{ and } \llVert U \rrVert_{2,0} \leq
R_0 \bigr\}, &\quad if $q = 0$.}
\]
\end{definition*}
The constraints on $R_q$ arise from the fact that the vector of row-wise
$\ell_2$ norms of a $p \times d$ orthonormal matrix belongs to a
sphere of radius
$d$. Roughly speaking, row sparsity asserts that there is a small
subset of variables
(coordinates of $\Real^p$)\vadjust{\goodbreak} that generate the principal subspace. Since
$\llVert\cdot\rrVert_{2,q}$ is coordinate-independent,
\emph{every}
orthonormal basis of a
$\mathcal{G}\in\mathcal{M}_{q} (R_q)$ has the same $(2,q)$-norm.

Another related notion of subspace sparsity is \emph{column sparsity},
which asserts that there is \emph{some} orthonormal basis of sparse vectors
that spans the principal subspace. Define the $(*,q)$-norm,
$q \in[0,\infty]$, of a $p \times d$ matrix $A$ as the maximal $\ell
_q$ norm
of its columns:
\[
\llVert A \rrVert_{*,q}:= \max_{1 \leq j \leq d} \llVert
a_{*j} \rrVert_q,
\]
where $a_{*j}$ denotes the $j$th column of $A$.
This is not coordinate-independent. We define the column sparse
subspaces to be those
that have some orthonormal basis with small $(*,q)$-norm.
\begin{definition*}[(Column sparse subspaces)]
For $0 \leq q < 2$ and $1 \leq R_q \leq p^{1-q/2}$,
\[
\mathcal{M}_{q}^{*} (R_q):= \cases{ \bigl\{
\colspan(U)\dvtx U \in\mathbb{V}_{p,d} \mbox{ and } \llVert U \rrVert
_{*,q}^q \leq R_q\bigr\}, &\quad if $0 < q < 2$\quad
and
\vspace*{2pt}\cr
\bigl\{ \colspan(U)\dvtx U \in\mathbb{V}_{p,d} \mbox{ and } \llVert
U \rrVert_{*,0} \leq R_0\bigr\}, &\quad if $q = 0$.}
\]
\end{definition*}
The column sparse subspaces are the $d$-dimensional subspaces that have
some orthonormal basis whose vectors are $\ell_q$ sparse in the usual sense.
Unlike row sparsity, the orthonormal bases of a column sparse
$\mathcal{G}$ do not all have the same $(*,q)$-norm, but if
$\mathcal{G} \in\mathcal{M}_{q}^{*} (R_q)$,
then there exists some $U \in\mathbb{V}_{p,d}$ such that
$\mathcal{G} = \colspan(U)$ and $\llVert U \rrVert_{*,q}^q
\leq R_q$
(or $\llVert U \rrVert_{*,q} \leq R_q$ for $q = 0$).

\subsection{Parameter space}
We assume that there exist i.i.d. random
vectors $Z_1,\ldots, Z_n \in\Real^p$, with $\E Z_1 = 0$ and $\Var
(Z_1) =
I_{p}$, such that
%
%
\begin{equation}
\label{eqsubgaussian-condition} X_i = \mu+ \Sigma^{1/2}
Z_i \quad\mbox{and}\quad\llVert Z_i \rrVert_{\psi_2}
\leq1
\end{equation}
for $i = 1,\ldots, n$, where
$\llVert\cdot\rrVert_{\psi_{\alpha}}$ is the Orlicz
$\psi_\alpha$-norm [e.g., \citet{vanderVaartAndWellner}, Chapter~2]
defined for $\alpha\geq1$ as
\[
\llVert Z \rrVert_{\psi_\alpha}:= \sup_{b: \llVert b \rrVert_2 \leq
1} \inf\biggl
\{C>0\dvtx\E\exp\biggl\llvert{ \frac{\langle Z,b\rangle
}{C}}\biggr\rrvert^{\alpha}\le2
\biggr\}.
\]
This ensures that all one-dimensional marginals of $X_i$ have
sub-Gaussian tails.
We also assume that the eigengap $\lambda_d - \lambda_{d+1} > 0$ so
that the principal subspace $\mathcal{S}$ is well defined.
Intuitively,
$\mathcal{S}$ is harder to estimate when the eigengap is
small. This
is made precise by the \emph{effective noise variance}
%
%
\begin{equation}
\label{eqnoise-to-signal} \sigma_d^2(
\lambda_1,\ldots,\lambda_p):= \frac{\lambda_1 \lambda_{d+1}}{(\lambda_d
- \lambda_{d+1})^2}.
\end{equation}
It turns out that this is a key quantity in the estimation of
$\mathcal{S}$, and that it is analogous to the noise
variance in
linear regression.
Let
\[
\mathcal{P}_{q} \bigl(\sigma^2,R_q\bigr)
\]
denote the class of
distributions on $X_1,\ldots,X_n$ that satisfy (\ref
{eqsubgaussian-condition}),
$\sigma_d^2 \leq\sigma^2$, and
$\mathcal{S}\in\mathcal{M}_{q} (R_q)$.
Similarly, let
\[
\mathcal{P}_{q}^* \bigl(\sigma^2, R_q\bigr)
\]
denote the class of distributions
that satisfy (\ref{eqsubgaussian-condition}),
$\sigma_d^2 \leq\sigma^2$, and
$\mathcal{S}\in\mathcal{M}_{q}^{*} (R_q)$.


\subsection{Subspace distance}
A notion of distance between subspaces is necessary to measure the performance
of a principal subspace estimator. The \emph{canonical
angles} between subspaces generalize the notion of angles between lines
and can be used to define subspace distances. There are several equivalent
ways to describe canonical angles, but for our purposes it will be
easiest to
describe them in terms of projection matrices. See
Bhatia [(\citeyear{Bhatia}), Chapter VII.1] and \citet{StewartAndSun}
for additional
background on canonical angles. For a subspace
$\mathcal{E} \in\mathbb{G}_{p,d}$ and its orthogonal projector $E$,
we write
$E^\perp$ to denote the orthogonal projector of $\mathcal{E}^\perp$ and
recall that $E^\perp= I_{p} - E$.
\begin{definition*}
Let $\mathcal{E}$ and $\mathcal{F}$ be $d$-dimensional subspaces of
$\Real^p$
with orthogonal projectors $E$ and $F$.
Denote the singular values of $E F^\perp$ by
$s_1 \geq s_2 \geq\cdots\,$. The \emph{canonical angles} between
$\mathcal{E}$ and
$\mathcal{F}$ are the numbers
\[
\theta_k(\mathcal{E}, \mathcal{F}) = \arcsin(s_k)
\]
for $k = 1,\ldots, d$ and the \emph{angle operator} between $\mathcal
{E}$ and
$\mathcal{F}$ is the $d \times d$ matrix
\[
\Theta(\mathcal{E}, \mathcal{F}) = \diag(\theta_1,\ldots,
\theta_d).
\]
\end{definition*}
In this paper we will consider the distance between subspaces $\mathcal
{E}, \mathcal{F} \in\mathbb{G}_{p,d}$
\[
\bigl\llVert\sin\Theta(\mathcal{E},\mathcal{F}) \bigr\rrVert_F,
\]
where $\llVert\cdot\rrVert_F$ is the Frobenius norm. This
distance is indeed a
metric on $\mathbb{G}_{p,d}$ [see \citet{StewartAndSun}, e.g.],
and can
be connected to the familiar Frobenius (squared error) distance between
projection matrices by the following fact from matrix
perturbation theory.
%
%
\begin{proposition}[{[See \citet{StewartAndSun}, Theorem I.5.5]}]
\label{procanonical-angles}
Let $\mathcal{E}$ and $\mathcal{F}$ be $d$-dimensional subspaces of
$\Real^p$
with orthogonal projectors $E$ and $F$.
Then:
\begin{enumerate}
\item The singular values of $E F^{\perp}$ are
\[
s_1,s_2,\ldots,s_d, 0,\ldots, 0.
\]
\item The singular values of $E - F$ are
\[
s_1,s_1, s_2, s_2,\ldots,s_d, s_d, 0,\ldots, 0.
\]
\end{enumerate}
In other words,
$E F^\perp$ has at most $d$ nonzero singular values
and the nonzero singular values of $E - F$ are
the nonzero singular values of $E F^{\perp}$,
each counted twice.\vadjust{\goodbreak}
\end{proposition}
Thus,
%
%
\begin{equation}
\label{eqsubspace-distance-identity} \bigl\llVert\sin\Theta(\mathcal{E},
\mathcal{F}) \bigr\rrVert_F^2 = \bigl\llVert E
F^{\perp} \bigr\rrVert_F^2 = \tfrac{1}{2}
\llVert E - F \rrVert_F^2 = \bigl\llVert
E^{\perp} F \bigr\rrVert_F^2.
\end{equation}
We will frequently use these identities.
For simplicity, we will overload notation and write
\[
\sin(U_1,U_2):=\sin\Theta\bigl( \colspan(U_1),
\colspan(U_2) \bigr)
\]
for $U_1,U_2 \in\mathbb{V}_{p,d}$. We also use a similar convention for
$\sin(E,F)$,
where $E, F$ are the orthogonal projectors corresponding to $\mathcal
{E},\mathcal{F} \in\mathbb{G}_{p,d}$
The following proposition, proved in Appendix \ref{secauxproofs},
relates the subspace
distance to the ordinary Euclidean distance between orthonormal matrices.
%
%
\begin{proposition}
\label{procosine-trace}
If $V_1,V_2 \in\mathbb{V}_{p,d}$, then
\[
\frac{1}{2} \inf_{Q \in\mathbb{V}_{d,d}} \llVert V_1 -
V_2 Q \rrVert_F^2 \leq\bigl\llVert
\sin(V_1, V_2) \bigr\rrVert_F^2
\leq\inf_{Q \in\mathbb{V}_{d,d}} \llVert V_1 - V_2 Q
\rrVert_F^2.
\]
\end{proposition}
In other words, the distance between two subspaces is equivalent to the
minimal distance between their orthonormal bases.

\subsection{Sparse subspace estimators}
\label{secsparse-subspace-estimators}
Here we introduce an estimator that achieves the optimal (up to a constant
factor) minimax error for row sparse subspace estimation. To estimate a
row sparse
subspace, it is natural to consider the empirical minimization problem
corresponding to (\ref{eqsubspace-pca}) with an additional sparsity
constraint corresponding to
$\mathcal{M}_{q} (R_q)$.

We define the row sparse principal subspace estimator to be a solution
of the
following constrained optimization problem:
%
%
\begin{eqnarray}
\label{eqrowsparse-subspace-pca-grassmann} &&\mbox{minimize}\quad \frac{1}{n}
\sum_{i=1}^n \bigl\llVert(I_{p}
- \Pi_{\mathcal{G}}) (X_{i} - \bar{X}) \bigr\rrVert
_2^2
\nonumber\\[-8pt]\\[-8pt]
&&\mbox{subject to}\quad \mathcal{G} \in\mathcal{M}_{q} (R_q).
\nonumber
\end{eqnarray}
For our analysis, it is more convenient to work on the Stiefel manifold.
Let $\langle A,B \rangle:=\tr(A^T B)$ for matrices $A,B$ of
compatible dimension. It is straightforward to show that following
optimization problem
is equivalent to (\ref{eqrowsparse-subspace-pca-grassmann}):
%
%
\begin{eqnarray}
\label{eqrowsparse-subspace-pca} &&\mbox{maximize}\quad \bigl\langle
\SampleCov,U
U^T \bigr\rangle
\nonumber
\\
&&\mbox{subject to}\quad U \in\mathbb{V}_{p,d}
\\
&&\hphantom{\mbox{subject to}\quad}\llVert U \rrVert_{2,q}^q \leq R_q
\qquad\mbox{(or $\llVert U \rrVert_{2,0} \leq R_0$ if $q=0$)}.
\nonumber
\end{eqnarray}
If $\hat{V}$ is a global maximizer of (\ref
{eqrowsparse-subspace-pca}), then
$\colspan(\hat{V})$ is a solution of
(\ref{eqrowsparse-subspace-pca-grassmann}). When $q = 1$, the estimator
defined by (\ref{eqrowsparse-subspace-pca}) is essentially a
generalization to
subspaces of the Lasso-type sparse PCA estimator proposed by
\citet{Jolliffe2003}.\vadjust{\goodbreak} A similar idea has also been used by
\citet
{Chen2010} in
the context of sufficient dimension reduction. 
The constraint set in (\ref{eqrowsparse-subspace-pca}) is clearly nonconvex,
however this is unimportant, because the objective function is convex
and we know that the maximum of a convex function over a set $D$
is unaltered if we replace $D$ by its convex hull. Thus,
(\ref{eqrowsparse-subspace-pca}) is equivalent to a convex \emph
{maximization} problem.
Finding a global maximum of convex maximization problems is
computationally challenging
and efficient algorithms remain to be developed.
Nevertheless, in the most popular case $q=1$, some algorithms have been
proposed with promising empirical performance
[\citet{Shen2008,Witten2009}].

We define the column sparse principal subspace
estimator analogously to the row sparse principal subspace estimator, using
the column sparse subspaces $\mathcal{M}_{q}^{*} (R_q)$ instead of
the row
sparse ones. This leads to the following equivalent Grassmann and Stiefel
manifold optimization problems:
%
%
\begin{eqnarray}
\label{eqcolsparse-subspace-pca-grassmann} &&\mbox{minimize}\quad \frac{1}{n}
\sum_{i=1}^n \bigl\llVert(I_{p}
- \Pi_{\mathcal{G}}) (X_{i} - \bar{X}) \bigr\rrVert
_2^2
\nonumber\\[-8pt]\\[-8pt]
&&\mbox{subject to}\quad \mathcal{G} \in\mathcal{M}_{q}^{*}
(R_q)
\nonumber
\end{eqnarray}
and
%
%
\begin{eqnarray}
\label{eqcolsparse-subspace-pca} &&\mbox{maximize}\quad \bigl\langle
\SampleCov,U
U^T \bigr\rangle
\nonumber
\\
&&\mbox{subject to}\quad U \in\mathbb{V}_{p,d}
\\
&&\hphantom{\mbox{subject to}\quad}\llVert U \rrVert_{*,q}^q \leq R_q
\qquad\mbox{(or $\llVert U \rrVert_{*,0} \leq R_0$ if $q=0$)}
\nonumber
\end{eqnarray}

\section{Main results}
\label{secresults}

In this section, we present our main results on the minimax lower and
upper bounds on sparse principal subspace estimation over the row
sparse and
column sparse classes.

\subsection{Row sparse lower bound}
To highlight the key results with minimal assumptions, we will
first consider the simplest case where $q=0$.
Consider the following two conditions.
%
%
\begin{condition}
\label{conconsistent-regime}
There is a constant $M > 0$ such that
\[
(R_q - d) \biggl[ \frac{\sigma^2}{n} \biggl( d + \log
\frac{(p-d)^{1-q/2}}{R_q-d} \biggr) \biggr]^{1 - {q/2}}
\leq M.
\]
\end{condition}
%
%
\begin{condition}
\label{connontrivial-rowsparse-parameters}
$4 \leq p - d$ and
$2d \leq R_q - d \leq(p - d)^{1 -{q/2}}$.
\end{condition}
Condition \ref{conconsistent-regime} is necessary for the existence
of a consistent estimator (see Theorems \ref
{lemlowerboundvarselecrow} and \ref{lemlowerboundparrow}).
Without Condition \ref{conconsistent-regime}, the statements of our
results would be complicated by multiple\vspace*{1pt} cases to deal with the fact
that the subspace distance is bounded above by $\sqrt{d}$.
The lower bounds on $p - d$ and $R_q - d$ are minor technical
conditions that
ensure our nonasymptotic bounds are nontrivial. Similarly, the upper
bound on
$R_q - d$ is only violated in trivial cases (detailed discussion given below).

%
%
\begin{theorem}[(Row sparse lower bound, $q = 0$)]
\label{thmlowerboundrow0}
If Conditions \ref{conconsistent-regime} and \ref
{connontrivial-rowsparse-parameters} hold, then
\[
\inf_{\hat{\mathcal{S}}} \sup_{\mathcal{P}_{0} (\sigma^2,R_0)} \E\bigl
\llVert\sin
\Theta(\hat{\mathcal{S}}, \mathcal{S}) \bigr\rrVert_F^2
\geq c (R_0 - d) \frac{\sigma^2}{n} \biggl[ d + \log\frac{p-d}{R_0-d}
\biggr].
\]
\end{theorem}
Here, as well as in the entire paper, $c$ denotes a universal, positive
constant, not necessarily the same at each occurrence.
This lower bound result reflects two separate aspects of the
estimation problem: \emph{variable selection}
and \emph{parameter estimation after variable selection}.
Variable selection refers to finding the variables that generate the principal
subspace, while estimation refers to estimating the subspace after selecting
the variables. For each variable, we accumulate two types of errors: one
proportional to $d$ that reflects the coordinates of the variable in the
$d$-dimensional subspace, and one proportional to $\log[(p-d)/(R_0 -
d)]$ that
reflects the cost of searching for the $R_0$ active variables. We prove
Theorem \ref{thmlowerboundrow0} in Appendix \ref{secprooflowerbound}.

The nonasymptotic lower bound for $0<q<2$ has
a more complicated dependence on ($n$, $p$, $d$, $R_q$, $\sigma^2$)
because of the interaction between $\ell_q$ and $\ell_2$ norms.
Therefore, our main lower bound result for $0<q<2$ will focus
on combinations of ($n$, $p$, $d$, $R_q$, $\sigma^2$) that correspond
to the
high-dimensional and sparse regime. (We state more general lower
bound results in Appendix~\ref{secprooflowerbound}.)
Let
%
%
\begin{equation}
\label{eqgamma-T-definition} T:=\frac{R_q - d}{(p-d)^{1-q/2}}
\quad\mbox{and}\quad\gamma:=
\frac{(p-d)\sigma^2}{n}.
\end{equation}
The interpretation for these two quantities is natural.
First, $T$ measures the relative sparsity of the problem. Roughly
speaking, it ranges between
$0$ and $1$ when the sparsity constraint in (\ref
{eqrowsparse-subspace-pca}) is active,
though the ``sparse'' regime generally corresponds to $T \ll1$.
The second quantity, $\gamma$ corresponds to the classic mean squared error
(MSE) of standard PCA.
The problem is low-dimensional if $\gamma$ is small compared to $T$.
We impose the following condition to preclude
this case.
%
%
\begin{condition}
\label{conrow-sparsity}
There is a constant $a < 1$ such that
$T^a \leq\gamma^{q/2}$.
\end{condition}
This condition lower bounds the classic MSE in terms of the sparsity
and is
mild in high-dimensional situations. When $a=q/2$, for example,
Condition \ref{conrow-sparsity} reduces to
\[
R_q-d \leq\frac{\sigma^2}{n} (p-d)^{2 - {q/2}}.
\]
We also note that this assumption becomes milder for larger values of $a$
and it is related to conditions in other minimax inference problem involving
$\ell_p$ and $\ell_q$ balls [see \citet{Donoho1994}, e.g.].

%
%
\begin{theorem}[(Row sparse lower bound, $0 < q < 2$)]
\label{thmlowerboundrowgeneral}
Let $q \in(0,2)$.
If Conditions \ref{conconsistent-regime} to
\ref{conrow-sparsity}
hold, then
\[
\inf_{\hat{\mathcal{S}}} \sup_{ \mathcal{P}_{q} (\sigma^2,R_q) } \E\bigl
\llVert\sin
\Theta(\hat{\mathcal{S}}, \mathcal{S}) \bigr\rrVert_F^2
\geq c (R_q-d) \biggl\{ \frac{\sigma^2}{n} \biggl[ d + \log
\frac{(p-d)^{1-q/2}}{R_q-d} \biggr] \biggr\}^{1 -{q/2}}.
\]
\end{theorem}
This result generalizes Theorem \ref{thmlowerboundrow0} and reflects
the same combination of variable selection and parameter estimation.
When Condition \ref{conrow-sparsity} does not hold, the problem is
outside of the
sparse, high-dimensional regime. As we show in the proof, there is
actually a
``phase transition regime'' between the high-dimensional sparse and the
classic dense regimes for which sharp minimax rate remains unknown. A similar
phenomenon has been observed in \citet{Birnbaum2012}.

\subsection{Row sparse upper bound}
Our upper bound results are obtained by analyzing the estimators given in
Section~\ref{secsparse-subspace-estimators}. The case where $q=0$ is
the clearest,
and we begin by stating a weaker, but simpler minimax upper bound for
the row
sparse class.
%
%
\begin{theorem}[(Row sparse upper bound, $q=0$)]
\label{thmupperboundrow0}
Let $\hat{\mathcal{S}}$ be any global maximizer of
(\ref{eqrowsparse-subspace-pca-grassmann}).
If $6\sqrt{R_0(d+\log p)} \leq\sqrt{n}$, then
\[
\sup_{\mathcal{P}_{0} (\sigma^2,R_0)} \E\bigl\llVert\sin\Theta(\hat
{\mathcal{S}},
\mathcal{S}) \bigr\rrVert_F^2 \leq c R_0
\frac{\lambda_1}{\lambda_{d+1}} \frac{\sigma^2(d+\log p)}{n}.
\]
\end{theorem}
Although (\ref{eqrowsparse-subspace-pca-grassmann}) may not have a
unique global optimum,
Theorem \ref{thmupperboundrow0} shows that \emph{any} global
optimum will be
within a certain radius of the principal subspace~$\mathcal{S}$.
The proof of
Theorem \ref{thmupperboundrow0},
given in Section~\ref{secquickupperboundproof}, is relatively
simple but still
nontrivial. It also serves as a prototype for the much more involved
proof of
our main upper bound result stated in
Theorem \ref{thmupperboundrowsparse} below. We note that
the rate given by Theorem~\ref{thmupperboundrow0} is off by a
$\lambda_1/\lambda_{d+1}$ factor that is due to the specific approach
taken to control an empirical process in our proof of
Theorem \ref{thmupperboundrow0}.

To state the main upper bound result with optimal dependence on
($n$, $p$, $d$, $R_q$, $\sigma^2$), we first describe some regularity
conditions.
Let
\[
\varepsilon_n:= \sqrt{2}R_q^{1/2} \biggl(
\frac{d+\log p}{n} \biggr)^{{1}/{2}
-{q}/{4}}.
\]
The regularity conditions are
%
%
\begin{eqnarray}
\label{eqsmallepsilon1} \varepsilon_n
&\leq&1,
\\
\label{eqsmallepsilon2}\quad c_1\sqrt{
\frac{d}{n}}\log n \lambda_1 + c_3
\varepsilon_n(\log n)^{5/2} \lambda_{d+1} &<&
\frac{1}{2}(\lambda_d-\lambda_{d+1}),
\\
\label{eqsmallepsilon3} c_3
\varepsilon_n(\log n)^{5/2}\lambda_{d+1} &\leq&\sqrt{
\lambda_1\lambda_{d+1}}^{1-q/2}(\lambda_d-
\lambda_{d+1})^{q/2}
\end{eqnarray}
and
\begin{equation}
\label{eqsmallepsilon4} c_3
\varepsilon_n^2(\log n)^{5/2}\lambda_{d+1}
\leq\sqrt{\lambda_1\lambda_{d+1}}^{2-q} (
\lambda_d-\lambda_{d+1})^{-(1-q)},
\end{equation}
where $c_1$ and $c_3$ are positive constants involved in the empirical
process arguments.
Equations
(\ref{eqsmallepsilon1}) to (\ref{eqsmallepsilon4}) require that
$\varepsilon_n$, the minimax rate of estimation (except the factor
involving $\lambda$),
to be small enough, compared to empirical process constants and
some polynomials of $\lambda$. Such conditions are mild in
the high dimensional, sparse regime, since to some extent, they
are qualitatively similar and analogous to Conditions \ref
{conconsistent-regime} to
\ref{conrow-sparsity} required by the lower bound.

%
%
\begin{remark}\label{remexplain-upper-bound-condition}
Conditions (\ref{eqsmallepsilon1}) to (\ref{eqsmallepsilon4})
are general enough to allow $R_q$, $d$ and $\lambda_j$ ($j=1,d,d+1$)
to scale
with $n$. For example, consider the case $q=0$, and let
$d=n^a$, $R_0=n^b$, $p=n^c$, $\lambda_1=n^{r_1}$, $\lambda_d=n^{r_2}$,
$\lambda_{d+1}=n^{r_3}$, where $0<a<b<c$, and $r_1\ge r_2 >r_3$.
Note that the $r_j$'s can be negative. Then it is straightforward to verify
that conditions (\ref{eqsmallepsilon1}) to (\ref{eqsmallepsilon4})
hold for large values of $n$ whenever $a+b <1$ and $r_1< r_2 + (1-a)/2$.
Condition (\ref{eqsmallepsilon1}) implies that $d$ cannot grow
faster than $\sqrt{n}$.
\end{remark}

%
%
\begin{theorem}[(Row sparse upper bound in probability)]
\label{thmupperboundrowsparse}
Let $q\in[0,1]$ and $\hat{\mathcal{S}}$ be any
solution of
(\ref{eqrowsparse-subspace-pca-grassmann}).
If $(X_1,\ldots,X_n) \sim\prob\in\mathcal{P}_{q} (\sigma^2,R_q)$
and
(\ref{eqsmallepsilon1}) to (\ref{eqsmallepsilon4})
hold, then
\[
\bigl\llVert\sin\Theta(\hat{\mathcal{S}}, \mathcal{S}) \bigr\rrVert
_F^2 \leq c R_q \biggl(\frac{\sigma^2(d+\log p)}{n}
\biggr)^{1-q/2}
\]
with probability at least $1-4/(n-1)-6\log n/n-p^{-1}$.
\end{theorem}

Theorem \ref{thmupperboundrowsparse} is presented in terms of a
probability bound
instead of an expectation bound. This stems from technical aspects of our
proof that involve bounding the supremum of an empirical process over a
set of
random diameter. For $q\in[0,1]$, the upper bound matches our lower bounds
(Theorems \ref{thmlowerboundrow0} and \ref{thmlowerboundrowgeneral})
for the entire tuple ($n$, $p$, $d$, $R_q$, $\sigma^2$) up to a constant
if
%
%
\begin{equation}
\label{eqadditional-condition} R_q^{2/(2-q)}\le
p^{c}
\end{equation}
for some constant $c<1$. To see this,
combining this additional condition and Condition
\ref{connontrivial-rowsparse-parameters}, the term $\log
\frac{p^{1-q/2}}{R_q}$ in the lower bound given in Theorem
\ref{thmlowerboundrowgeneral} is within a constant factor of $\log
p$ in the upper bound given in Theorem~\ref{thmupperboundrowsparse}.
It is straightforward to check that
the other terms in lower and upper bounds agree up to constants with
obvious correspondence. Moreover, we note that the additional condition
(\ref{eqadditional-condition}) is only slightly stronger than the last
inequality in Condition \ref{connontrivial-rowsparse-parameters}. The
proof of Theorem \ref{thmupperboundrowsparse} is in Appendix~\ref{secproofupperboundmain}.

Using the probability upper bound result and the fact that
$\llVert{\sin\Theta}(\hat{\mathcal{S}},
\mathcal{S}) \rrVert_F^2
\leq d$, one can derive an upper bound
in expectation.
%
%
\begin{corollary}
\label{corupperbound-expectation}
Under the same condition as in Theorem \ref
{thmupperboundrowsparse}, we have
for some constant $c$,
\[
\E\bigl\llVert\sin\Theta(\hat{\mathcal{S}}, \mathcal{S}) \bigr\rrVert
_F^2 \leq c \biggl\{ R_q \biggl[
\frac{\sigma^2(d+\log p)}{n} \biggr]^{1-q/2}+d \biggl( \frac
{\log n}{n}+
\frac{1}{p} \biggr) \biggr\}.
\]
\end{corollary}

%
%
\begin{remark}
\label{remmatching-expectation}
The expectation upper bound has an additional
$
d(\log n / n + 1/p)
$
term that can be further reduced by refining the argument (see
Remark \ref{remupperbound-expectation} below). It is not obvious
if one
can completely avoid such a term. But in many situations it
is dominated by the first term. Again, we invoke the scaling considered in
Remark \ref{remexplain-upper-bound-condition}. When $q=0$, the
first term
is of order $n^{a+b+(r_1+r_3)/2-r_2-1}$, and the additional term is
$n^{a-1}\log n + n^{a-c}$, which is asymptotically negligible if
$b>r_2-(r_1+r_3)/2+(1-c)_+$.
\end{remark}

%
%
\begin{remark}
\label{remupperbound-expectation}
Given any $r>0$,
it is easy to modify the proof of Theorem \ref
{thmupperboundrowsparse} [as well
as conditions (\ref{eqsmallepsilon1}) to (\ref{eqsmallepsilon4})]
such that the results of Theorem \ref{thmupperboundrowsparse} and
Corollary~\ref{corupperbound-expectation} hold with $c$ replaced
by some constant $c(r)$,
and the probability bound becomes $1-4/(n^r-1)-6\log n/ n^{r} - 1/p^r$.
\end{remark}
%

\subsection{Column sparse lower bound}

By modifying the proofs of
Theorems \ref{thmlowerboundrow0} and \ref
{thmlowerboundrowgeneral}, we can obtain lower bound
results for the column sparse case that are parallel to the row sparse case.
For brevity, we present the $q=0$ and $q>0$ cases together. The analog
of $T$,
the degree of sparsity,
for the column sparse case is
%
%
\begin{equation}
\label{eqTstar-definition} T_*:=\frac{d(R_q - 1)}{(p-d)^{1-q/2}},
\end{equation}
and the analogs of
Conditions \ref{connontrivial-rowsparse-parameters} and \ref
{conrow-sparsity} are the following.
%
%
\begin{condition}
\label{connontrivial-columnsparse-parameters}
$4 d \leq p - d$ and
$d \leq d(R_q - 1) \leq(p-d)^{1 -{q/2}}$.
\end{condition}
%
%
\begin{condition}
\label{concolumn-sparsity}
There is a constant $a < 1$ such that $T_*^a \leq\gamma^{q/2}$.
\end{condition}

%
%
\begin{theorem}[(Column sparse lower bound)]
\label{thmlowerboundcol}
Let $q \in[0,2)$.
If Conditions \ref{connontrivial-columnsparse-parameters} and \ref
{concolumn-sparsity}
hold, then
\[
\inf_{\hat{\mathcal{S}}} \sup_{ \mathcal{P}_{q}^* (\sigma^2,R_q) } \E
\bigl\llVert\sin(
\hat{\mathcal{S}}, \mathcal{S}) \bigr\rrVert_F^2
\geq c d(R_q-1) \biggl\{ \frac{\sigma^2}{n} \biggl[ 1 + \log
\frac{(p-d)^{1-q/2}}{d(R_q-1)} \biggr] \biggr\}^{1 -{q/2}}.
\]
\end{theorem}
For column sparse subspaces, the lower bound is dominated by the variable
selection error, because column sparsity is defined in terms of the maximal
$\ell_0$ norms of the vectors in an orthonormal basis and $R_0$
variables must
be selected for each of the $d$ vectors. So the variable selection
error is
inflated by a factor of~$d$, and hence becomes the dominating term in
the total
estimation error. We prove Theorem \ref{thmlowerboundcol} in
Appendix \ref{secprooflowerbound}.

\subsection{Column sparse upper bound}

A specific challenge in analyzing the column sparse principal subspace
problem (\ref{eqcolsparse-subspace-pca})
is to bound the
supremum of the empirical process
\[
\bigl\langle\SampleCov-\Sigma,UU^T-VV^T\bigr\rangle
\]
indexed by all $U\in\mathcal U(p,d,R_q,\varepsilon)$ where
\[
\mathcal U(p,d,R_q,\varepsilon)\equiv\bigl\{U\dvtx\mathbb{V}_{p,d},
\llVert U \rrVert_{*,q}^q\le R_q, \bigl\llVert
UU^T-VV^T \bigr\rrVert_F\le\varepsilon\bigr\}.
\]
Unlike the row sparse matrices, the matrices $UU^T$ and $VV^T$ are no longer
column sparse with the same radius $R_q$.

By observing that
$\mathcal{M}_{q}^{*} (R_q)\subseteq\mathcal{M}_{q} (dR_q)$,
we can reuse the proof of Theorem~\ref{thmupperboundrowsparse} to
derive the following
upper bound for the column sparse class.

%
\begin{corollary}[(Column sparse upper bound)]
\label{thmuppercolumn}
Let $q\in[0,1]$ and $\hat{\mathcal{S}}$ be any
solution of
(\ref{eqcolsparse-subspace-pca-grassmann}).
If $(X_1,\ldots,X_n) \sim\prob\in\mathcal{P}_{q}^* (\sigma^2,R_q)$
and
(\ref{eqsmallepsilon1}) to (\ref{eqsmallepsilon4})
hold with $R_q$ replaced by $d R_q$, then
\[
\bigl\llVert\sin\Theta(\hat{\mathcal{S}}, \mathcal{S}) \bigr\rrVert
_F^2 \leq c d R_q \biggl(\frac{\sigma^2(d+\log p)}{n}
\biggr)^{1-q/2}
\]
with probability at least $1-4/(n-1)-6\log n/n-p^{-1}$.
\end{corollary}
Corollary \ref{thmuppercolumn} is slightly weaker than the
corresponding result for
the row sparse class. It matches the lower bound in Theorem \ref
{thmlowerboundcol}
up to a constant if
\[
\bigl(d(R_q-1)\bigr)^{2/(2-q)} \leq p^c
\]
for some constant $c < 1$,
and $d < C \log p$ for some other constant $C$.

\subsection{A conjecture for the column sparse case}
Note that Theorem \ref{thmlowerboundcol} and Corollary \ref
{thmuppercolumn} only match
when $d\le C\log p$. For larger values of $d$, we believe that
the lower bound in Theorem \ref{thmlowerboundcol} is optimal and
the upper bound
can be improved.\vadjust{\goodbreak}

\begin{Conjecture*}[(Minimax error bound for column sparse
case)] Under the same conditions
as in Corollary \ref{thmuppercolumn}, there exists an estimator
$\hat{\mathcal{S}}$ such
that
\[
\bigl\llVert\sin\Theta(\hat{\mathcal{S}}, \mathcal{S}) \bigr\rrVert
_F^2 \leq c d R_q \biggl(\frac{\sigma^2(1+\log p)}{n}
\biggr)^{1-q/2}
\]
with high probability. As a result, the optimal minimax lower and upper bounds
for this case shall be
\[
\bigl\llVert\sin\Theta(\hat{\mathcal{S}}, \mathcal{S}) \bigr\rrVert
_F^2\asymp d R_q \biggl(\frac{\sigma^2\log p}{n}
\biggr)^{1-q/2}.
\]
\end{Conjecture*}

One reason for the conjecture is based on the following intuition.
Suppose that
$\lambda_1>\lambda_2>\cdots>\lambda_d>\lambda_{d+1}$ (there is enough
gap between the leading
eigenvalues) one can recover
the individual leading eigenvectors with an error rate whose dependence
on $(n,R_q,p)$ is the same
as in the lower bound [cf. \citet{Vu2012,Birnbaum2012}].
As a result, the estimator $\hat V=(\hat v_1,\hat v_2,\ldots,\hat v_d)$
shall give the
desired upper bound. On the other hand, it remains open to us whether
the estimator
in (\ref{eqcolsparse-subspace-pca}) can achieve this rate for $d$
much larger than $\log p$.


\section{Sketch of proofs}
\label{secsketchproof}

For simplicity, we focus on the row sparse case with $q=0$,
assuming also the high dimensional and sparse regime.
For more general cases, see Theorems \ref
{lemlowerboundvarselecrow} and \ref{lemlowerboundparrow} in
Appendix \ref{secprooflowerbound}.

\subsection{The lower bound}\label{subsecsketch-lower}
Our proof of the lower bound features a combination of the general
framework of the Fano method
and a careful combinatorial analysis of packing sets of various classes
of sparse matrices.
The particular challenge is to construct a rich packing set of
the parameter space $\mathcal{P}_{q} (\sigma^2,R_q)$. We will
consider centered $p$-dimensional Gaussian distributions with
covariance matrix $\Sigma$ given by
%
%
\begin{equation}
\label{eqspiked-covariance-model-1} \Sigma(A)=b AA^T +
I_p,
\end{equation}
where
$A\in\mathbb{V}_{p,d}$ is constructed from the ``local Stiefel
embedding'' as given
below.
Let $1 \leq k \leq d < p$ and the function $A_\varepsilon\dvtx\mathbb
{V}_{p-d,k}
\mapsto\mathbb{V}_{p,d}$ be defined in block form as
%
%
\begin{equation}
\label{eqstiefel-embedding-1} A_\varepsilon(J) = \left[\matrix{ \bigl(1-
\varepsilon^2\bigr)^{1/2} I_{k} & 0
\cr
0 &
I_{d-k}
\cr
\varepsilon J & 0}\right]
\end{equation}
for $0 \leq\varepsilon\leq1$.
We have the following generic method for
lower bounding the minimax risk of estimating the principal subspace of a
covariance matrix. It is proved in Appendix \ref{secprooflowerbound}
as a consequence
of Lemmas \ref{lemfano} to \ref{lemstiefel-embedding}.
%
%
\begin{lemma}[(Fano method with Stiefel embedding)]
\label{lemfano-generic-stiefel-lower-bound}
Let $\varepsilon\in[0, 1]$ and $\{J_1,\ldots,J_N\} \subseteq\mathbb
{V}_{p-d,k}$
for $1 \leq k \leq d < p$.
For each $i = 1,\ldots, N$,
let $\prob_i$ be the $n$-fold product of the $\normal(0, \Sigma
(A_\varepsilon(J_i)))$
probability measure,
where
$\Sigma(\cdot)$ is defined in (\ref{eqspiked-covariance-model-1})
and
$A_\varepsilon(\cdot)$ is defined in (\ref{eqstiefel-embedding-1}).
If
\[
\min_{i \neq j} \llVert J_i - J_j
\rrVert_F \geq\delta_N,
\]
then every estimator $\hat{\mathcal{A}}$ of
$\mathcal{A}_i:=\colspan(A_\varepsilon(J_i))$ satisfies
\[
\max_i \E_i \bigl\llVert\sin\Theta(\hat{
\mathcal{A}}, \mathcal{A}_i) \bigr\rrVert_F \geq
\frac{\delta_N \varepsilon\sqrt{1-\varepsilon^2}}{2} \biggl[ 1 - \frac{
4 n k\varepsilon^2 / \sigma^2 + \log2}{\log N} \biggr],
\]
where $\sigma^2 = (1+b) / b^2$.
\end{lemma}
Note that if $\llVert J \rrVert_{2,0} \leq R_0-d$, then
$\llVert A_\varepsilon(J) \rrVert_{2,0} \leq R$.
Thus Lemma \ref{lemfano-generic-stiefel-lower-bound}
with appropriate choices of $J_i$ can yield minimax lower
bounds over $p$-dimensional Gaussian distributions whose
principal subspace is $R_0$ row sparse.

The remainder of the proof consists of two applications of
Lemma \ref{lemfano-generic-stiefel-lower-bound} that correspond to the
two terms in Theorem \ref{thmlowerboundrow0}.
In the first part, we use a variation of the Gilbert--Varshamov bound
(Lemma \ref{lemhypercube}) to construct a packing set in $\mathbb{V}_{p-d,1}$
consisting of $(R_0-d)$-sparse vectors.
Then we apply Lemma \ref{lemfano-generic-stiefel-lower-bound} with
\[
k = 1,\qquad \delta_N = 1/4,\qquad \varepsilon^2\asymp
\frac{\sigma^2 R_0\log p}{n}.
\]
This yields a minimax lower bound
that reflects the variable selection complexity.
In the second part, we leverage existing results on the metric
entropy of the Grassmann manifold (Lemma \ref{lemgrassmann-packing})
to construct a packing set of $\mathbb{V}_{R_0-d,d}$.
Then we apply Lemma \ref{lemfano-generic-stiefel-lower-bound}
with
\[
k = d,\qquad \delta_N = c_0\sqrt{d}/e,\qquad \varepsilon^2
\asymp\frac{\sigma^2 R_0}{n}.
\]
This yields a minimax lower bound
that reflects the complexity of
post-selection estimation. Putting these two results together,
we have for a subset of Gaussian distributions
$G \subseteq\mathcal{P}_{0} (\sigma^2,R_q)$
the minimax lower bound:
\[
\max_{G} \E\bigl\llVert\sin\Theta(\hat{\mathcal{A}},
\mathcal{A}_i) \bigr\rrVert_F^2 \geq c
R_0 \frac{\sigma^2}{n} (\log p \land d) \geq(c/2) R_0
\frac{\sigma^2}{n} (d + \log p).
\]

\subsection{The upper bound}\label{subsecsketch-upper}\label
{secquickupperboundproof}
The upper bound proof requires a careful analysis of the
behavior of the empirical maximizer of the PCA problem
under sparsity constraints.
The first key ingredient is to provide a lower bound of the
curvature of the objective function at its global maxima.
Traditional results of this kind, such as Davis--Kahan sin$\Theta$ theorem
and Weyl's inequality, are not sufficient for our purpose.

The following lemma, despite its elementary form, has not been seen in
the literature (to our knowledge). It gives us the right tool to bound
the curvature of
the matrix functional $F\mapsto\langle A,F\rangle$ at its point of maximum
on the Grassmann manifold.

%
%
\begin{lemma}[(Curvature lemma)]
\label{lemcurvature-lemma}
Let $A$ be a $p\times p$ positive semidefinite matrix and
suppose that its eigenvalues $\lambda_1(A) \geq\cdots\geq\lambda_p(A)$
satisfy $\lambda_d(A) > \lambda_{d+1}(A)$ for $d < p$.
Let $\mathcal{E}$ be the $d$-dimensional subspace spanned by the
eigenvectors of $A$ corresponding to its $d$ largest eigenvalues, and let
$E$ denote its orthogonal projector. If $\mathcal{F}$
is a $d$-dimensional subspace of $\Real^p$ and $F$ is its orthogonal projector,
then
\[
\bigl\llVert\sin\Theta(\mathcal{E}, \mathcal{F}) \bigr\rrVert
_F^2 \leq\frac{ \langle A,E - F\rangle}{ \lambda_d(A) - \lambda
_{d+1}(A) }.
\]
\end{lemma}
Lemma \ref{lemcurvature-lemma} is proved in Appendix \ref
{subsecadditionalproofupper}.
An immediate corollary is the following alternative to the traditional matrix
perturbation approach to bounding subspace distances using the Davis--Kahan
$\sin\Theta$ theorem and Weyl's inequality.
%
%
\begin{corollary}[(Variational $\sin\Theta$)]
\label{lemcurvature-centered}
In addition to the hypotheses of Lem\-ma~\ref{lemcurvature-lemma},
if $B$ is a symmetric matrix and $F$ satisfies
%
%
\begin{equation}
\label{eqmonotone-condition} \langle B,E\rangle- g(E) \leq\langle
B,F\rangle-
g(F)
\end{equation}
for some function $g\dvtx\Real^{p \times p} \mapsto\Real$,
then
%
%
\begin{equation}
\label{eqcurvature-centered-bound} \bigl\llVert\sin\Theta(\mathcal{E},
\mathcal{F}) \bigr\rrVert_F^2 \leq\frac
{\langle B - A,F - E\rangle- [g(F) - g(E)]} {
\lambda_d(A) - \lambda_{d+1}(A)}.
\end{equation}
\end{corollary}
The corollary is different from the Davis--Kahan $\sin\Theta$ theorem
because the
orthogonal projector $F$ does not have to correspond to a subspace
spanned by
eigenvectors of $B$. $F$ only has to satisfy
\[
\langle B,E\rangle- g(E) \leq\langle B,F\rangle- g(F).
\]
This condition is suited ideally for analyzing solutions of regularized and/or
constrained maximization problems where $E$ and $F$ are feasible, but
$F$ is
optimal. In the simplest case, where $g \equiv0$, combining
(\ref{eqcurvature-centered-bound}) with the Cauchy--Schwarz inequality and
(\ref{eqsubspace-distance-identity}) recovers a form of the
Davis--Kahan $\sin\Theta$ theorem in the Frobenius norm:
\[
\frac{1}{\sqrt{2}} \bigl\llVert\sin\Theta(\mathcal{E}, \mathcal{F})
\bigr
\rrVert_F \leq\frac
{\llVert B-A \rrVert_F} {
\lambda_d(A) - \lambda_{d+1}(A)}.
\]

In the upper bound proof, let $V\in\mathbb{V}_{p,d}$ be the true parameter,
and $\hat V$ be a solution of (\ref{eqrowsparse-subspace-pca}).
Then we have
\[
\bigl\langle S_n,\hat V\hat V^T-V V^T\bigr
\rangle\ge0.
\]

Applying Corollary \ref{lemcurvature-centered} with
$B=\SampleCov$, $A=\Sigma$, $E=VV^T$, $F=\hat V\hat V^T$, and
$g\equiv0$,
we have
%
%
\begin{equation}
\label{equpperboundcurvature} \bigl\llVert\sin
\Theta(V, \hat V) \bigr\rrVert_F^2 \le
\frac{ \langle\SampleCov-\Sigma,\hat V\hat V^T - V V^T \rangle}{
\lambda_d(\Sigma) - \lambda_{d+1}(\Sigma) }.
\end{equation}

Obtaining a sharp upper bound for $\langle S-\Sigma,\hat V\hat V^T - V
V^T \rangle$ is nontrivial. First, one needs to control
$\sup_{F\in\mathcal F}\langle S-\Sigma,F \rangle$ for some class
$\mathcal F$ of sparse and symmetric matrices. This requires some
results on quadratic form empirical process. Second, in order to obtain
better bounds, we need to take advantage of the fact that $\hat V\hat
V^T-V V^T$ is probably small. Thus, we need to use a peeling argument
to deal with the case where $\mathcal F$ has a random (but probably)
small diameter. These details are given in Appendices
\ref{subsecproofmainupper} and \ref{secempiricalprocessproofs}. Here we
present a short proof of Theorem \ref{thmupperboundrow0} to illustrate
the idea.

\begin{pf*}{Proof of Theorem \ref{thmupperboundrow0}}
By (\ref{equpperboundcurvature}), we have
\[
\hat{\varepsilon}^2:= \bigl\llVert\sin\Theta(\hat{\mathcal{S}},
\mathcal{S}) \bigr\rrVert_F^2 \leq\frac
{\langle\SampleCov- \Sigma,\hat{V} \hat{V}^T - V V^T \rangle} {
\lambda_d - \lambda_{d+1}}
\]
and
%
%
\begin{equation}
\label{eqcurvature0} \hat\varepsilon^2 \leq
\frac{\sqrt{2}}{\lambda_d-\lambda_{d+1}} \biggl\langle{\SampleCov-\Sigma
} { \frac{ \hat V\hat V^T-VV^T} {
\llVert\hat V\hat V^T-VV^T \rrVert_F} } \biggr
\rangle\hat\varepsilon,
\end{equation}
because $\llVert\hat V\hat V^T-VV^T \rrVert_F^2 = 2 \hat
{\varepsilon}^2$ by
(\ref{eqsubspace-distance-identity}).
Let
\[
\Delta=\frac{ \hat V\hat V^T-VV^T}{\llVert\hat V\hat V^T-VV^T \rrVert_F}.
\]
Then $\llVert\Delta\rrVert_{2,0}\le2R_0$, $\llVert
\Delta\rrVert_F=1$,
and $\Delta$ has at most $d$ positive eigenvalues and
at most $d$ negative eigenvalues (see Proposition \ref{procanonical-angles}).
Therefore, we can write $\Delta=AA^T-BB^T$ where
$\llVert A \rrVert_{2,0}\le2R_0$, $\llVert A \rrVert_F\le1$, $A\in
\Real^{p\times d}$,
and the same holds for~$B$.
Let
\[
\mathcal{U}(R_0) = \bigl\{ U \in\Real^{p \times d}\dvtx\llVert U
\rrVert_{2,0} \leq2 R_0 \mbox{ and } \llVert U \rrVert
_F \leq1 \bigr\}.
\]
Equation (\ref{eqcurvature0}) implies
\[
\E\hat\varepsilon\leq\frac{2\sqrt{2}}{\lambda_d-\lambda_{d+1}} \E\sup
_{U\in\mathcal{U}(R_0)} \bigl
\llvert\bigl\langle\SampleCov-\Sigma,UU^T \bigr\rangle\bigr\rrvert.
\]
The empirical process $\langle\SampleCov-\Sigma,UU^T \rangle$
indexed by $U$ is a generalized quadratic form, and
a sharp bound of its supremum involves some
recent advances in empirical process theory due to \citet{Mendelson2010}
and extensions of his results.
By Corollary 4.1 of \citet{Vu2012b}, we have
\begin{eqnarray*}
&& \E\sup_{U \in\mathcal{U}(R_0)} \bigl\llvert\bigl\langle\SampleCov-
\Sigma,UU^T \bigr\rangle\bigr\rrvert
\\
&&\qquad\le c\lambda_1 \biggl\{ \frac{\E\sup_{U\in\mathcal
{U}(R_0)}\langle
\mathcal{Z},U \rangle
}{\sqrt{n}} + \biggl(
\frac{\E\sup_{U\in\mathcal{U}(R_0)}\langle\mathcal{Z},U \rangle
}{\sqrt{n}} \biggr)^2 \biggr\},
\end{eqnarray*}
where $\mathcal Z$ is a $p\times d$ matrix of i.i.d. standard Gaussian
variables.
To control $\E\sup_{U\in\mathcal{U}}\langle\mathcal{Z},U \rangle$,
note that
\[
\langle\mathcal{Z},U \rangle\le\llVert\mathcal Z \rrVert_{2,\infty}
\llVert U \rrVert_{2,1} \le\llVert\mathcal Z \rrVert_{2,\infty}
\sqrt{2R_0},
\]
because
$U\in\mathcal{U}(R_0)$. Using a standard
$\delta$-net argument (see Propositions~\ref{procovering-argument}
and~\ref{prosphere-covering}), we have, when
$p>5$,
%
%
\begin{equation}
\label{eqgaussian-row-norm-bound} \bigl\llVert{\llVert\mathcal{Z}
\rrVert
_{2,\infty}}\bigr\rrVert_{\psi
_2}\le4.15\sqrt{d+\log p}
\end{equation}
and hence
\[
\E\sup_{U\in\mathcal{U}} \langle\mathcal{Z},U \rangle\leq6
\sqrt{R_0(d+\log p)}.
\]
The proof is complete since we assume that $6\sqrt{R_0(d+\log p)}\le
\sqrt{n}$.
\end{pf*}

\section{Discussion}
\label{secdiscussion}

There is a natural correspondence between the sparse principal subspace
optimization
problem (\ref{eqrowsparse-subspace-pca-grassmann}) and some
optimization problems
considered in the sparse regression literature.
We have also found that there is a correspondence between minimax
results for sparse regression and those that we presented in this
article. In spite of these connections, results on computation for
sparse principal subspaces (and sparse PCA) are far less developed than
for sparse regression. In this final section, we will discuss the connections
with sparse regression, both optimization and minimax theory, and then
conclude with some open problems for sparse principal subspaces.

\subsection{Connections with sparse regression}
Letting $\tilde{X}_i = X_i - \bar{X}$ denote a centered observation,
we can write (\ref{eqrowsparse-subspace-pca-grassmann}) in the $d=1$
case as an equivalent penalized regression problem:
\begin{eqnarray*}
&&\mbox{minimize}\quad \frac{1}{n} \sum
_{i=1}^n \bigl\llVert\tilde{X}_i - u
u^T \tilde{X}_i \bigr\rrVert_2^2
+ \tau_q \llVert u \rrVert_q^q
\nonumber\\[-8pt]\\[-8pt]
&&\mbox{subject to}\quad u \in\Real^p \quad\mbox{and}\quad \llVert u \rrVert
_2=1\nonumber
\end{eqnarray*}
for $0 < q \leq1$ and similarly for $q=0$. The penalty parameter
$\tau_q \geq0$ plays a similar role as $R_q$. When $q=1$ this is
equivalent to a penalized form of the sparse PCA estimator considered
in \citet{Jolliffe2003} and it also bears similarity to the estimator
considered by \citet{Shen2008}. It is also similar to the famous
$\ell_1$-penalized optimization often used in high-dimensional
regression [\citet{Tibshirani1996}]. In the subspace case $d>1$, one
can write an analogous penalized multivariate regression problem:
\begin{eqnarray*}
&&\mbox{minimize}\quad \frac{1}{n} \sum
_{i=1}^n \bigl\llVert\tilde{X}_i - U
U^T \tilde{X}_i \bigr\rrVert_2^2
+ \tau_q \llVert U \rrVert_{2,q}^q
\\
&&\mbox{subject to}\quad U \in\mathbb{V}_{p,d}
\end{eqnarray*}
for $0 < q \leq1$ and similarly for $q=0$. When $q=1$, this
corresponds to a ``group Lasso'' penalty where entries in the same row
of $U$ are penalized simultaneously [\citet{Yuan2006,Zhao2009}]. The
idea being that as $\tau_q$ varies, a variable should enter/exit all
$d$ coordinates simultaneously. In the column sparse case, when $q=1$
the analogous penalized multivariate regression problem has a penalty
which encourages each column of $U$ to be sparse, but does not require
that the pattern of sparsity to be the same across columns.

The analogy between row sparse principal subspace estimation and sparse
regression goes beyond the optimization problems formulated above---it
is also reflected in terms of the minimax rate.
In the sparse regression problem, we assume an i.i.d. sample $(X_i,
Y_i)\in\Real^p\times\Real$ for $1\le i\le n$ satisfying
\[
Y_i=\beta^T X_i + \varepsilon_i,
\]
where $\varepsilon_i$ is mean zero, independently of $X_i$, and $\beta
\in\Real^p$
is the regression coefficient vector.
\citet{Raskutti2011} showed (with some additional conditions on the
distribution of $X_i$) that if $\llVert\beta\rrVert
_q^q\le R_q$ and
$\operatorname{Var}\varepsilon_i\le\sigma^2$, then
the minimax rate of estimating $\beta$ in $\ell_2$ norm is (ignoring
constants)
\[
\sqrt{R_q} \biggl(\frac{\sigma^2\log p}{n} \biggr)^{
{1}/{2}-{q}/{4}}.
\]
The estimator that achieves this error rate is obtained by solving the
$\ell_q$ constrained least square
problem.
The $d
> 1$ case corresponds to the multivariate regression model, where
$Y_i\in\Real^d$, $\beta\in
\Real^{p\times d}$, and
$\operatorname{Var}\varepsilon_i=\sigma^2I_d$.
The results of \citet{Negahban2012}, with straightforward modifications,
imply that if $\llVert\beta\rrVert_{2,q}^q\le R_q$, then
a penalized least squares estimator can achieve the $\ell_2$ error rate
\[
\sqrt{R_q} \biggl(\frac{\sigma^2(d+\log p)}{n} \biggr)^{
{1}/{2}-{q}/{4}},
\]
agreeing with our minimax lower and upper bounds for the row sparse
principal subspace problem.



\subsection{Practical concerns}
The nature of this work is theoretical and it leaves open many challenges
for methodology and practice. The minimax optimal estimators that we
present appear to be computationally intractable because they involve convex
\emph{maximization} rather than convex \emph{minimization} problems.
Even in
the case $q=1$, which corresponds to a subspace extension of $\ell_1$
constrained
PCA, the optimization problem remains challenging as there are no known
algorithms to efficiently compute a global maximum.

Although the
minimax optimal estimators that we propose do not require knowledge of the
noise-to-signal ratio $\sigma^2$, they do require knowledge of (or an
upper bound on) the sparsity $R_q$. It is not hard to modify our
techniques to
produce an estimator that gives up adaptivity to $\sigma^2$ in
exchange for
adaptivity to $R_q$. One could do this by using penalized versions of our
estimators with a penalty factor proportional to $\sigma^2$. An extension
along this line has already been considered by
\citet{Lounici2012} for the $d=1$ case. A more interesting question
is whether or not there exist fully adaptive principal subspace estimators.

Under what conditions can one find an estimator that achieves the
minimax optimal error without requiring knowledge of either $\sigma^2$
or $R_q$? Works by \citet{Paul2007} and \citet{Ma2011} on
refinements of diagonal thresholding for the spiked covariance model
seems promising on this front, but as we mentioned in the
\hyperref[secintroduction]{Introduction}, the spiked covariance model
is restrictive and necessarily excludes the common practice of
standardizing variables. Is it possible to be adaptive outside the
spiked covariance model? One possible approach can be described in the
following three steps. (1) use a conservative choice of $R_q$ (say,
$p^a$, for some $0<a<1$); (2) estimate $\sigma^2$ using eigenvalues
obtained from the sparsity constrained principal subspace estimator;
and (3) use a sparsity penalized principal subspace estimator with
$\sigma^2$ replaced by its estimate. We will pursue this idea in
further detail in future work.

%
%
\begin{appendix}\label{app}
\section{Lower bound proofs}
\label{secprooflowerbound}

Theorems \ref{thmlowerboundrow0}, \ref{thmlowerboundrowgeneral} and
\ref{thmlowerboundcol} are
consequences of three more general results stated below.
An essential part of the strategy of our proof is to analyze the \emph{variable
selection} and \emph{estimation} aspects of the problem separately.
We will consider two types of subsets of the parameter space that capture
the essential difficulty of each aspect: one where the subspaces vary over
different subsets of variables, and another where the subspaces vary
over a
fixed subset of variables. The first two results give lower bounds for
each aspect
in the row sparse case. Theorems \ref{thmlowerboundrow0} and
\ref{thmlowerboundrowgeneral}
follow easily from them. The third result directly addresses the proof of
Theorem \ref{thmlowerboundcol}.

%
%
\begin{theorem}[(Row sparse variable selection)]
\label{lemlowerboundvarselecrow}
Let $q \in[0,2)$ and $(p,d,R_q)$ satisfy
\[
4 \leq p - d \quad\mbox{and}\quad1 \leq R_q - d \leq(p - d)^{1-q/2}.
\]
There exists a universal constant $c > 0$ such that every estimator
$\hat{\mathcal{S}}$ satisfies the following.
If $T < \gamma^{q/2}$, then
%
%
\begin{eqnarray}
\label{eqsparse-variable-selection-lower-bound} && \sup_{\mathcal{P}_{q}
(\sigma^2,R_q)}
\E\bigl\llVert\sin\Theta(\hat{\mathcal{S}}, \mathcal{S}) \bigr\rrVert
_F
\nonumber\\[-8pt]\\[-8pt]
&&\qquad\geq c \biggl\{ (R_q - d) \biggl[ \frac{\sigma^2}{n} \bigl( 1 -
\log\bigl(T / \gamma^{q/2} \bigr) \bigr) \biggr]^{1-
{q/2}} \land1
\biggr\}^{1/2}.
\nonumber
\end{eqnarray}
Otherwise,
%
%
\begin{equation}
\label{eqdense-variable-selection-lower-bound} \sup_{\mathcal{P}_{q}
(\sigma^2,R_q)} \E
\bigl\llVert\sin\Theta(\hat{\mathcal{S}}, \mathcal{S}) \bigr\rrVert
_F \geq c \biggl\{ \frac{(p-d) \sigma^2}{n} \land1 \biggr
\}^{1/2}.
\end{equation}
\end{theorem}
The case $q = 0$ is
particularly simple, because $T < \gamma^{q/2} = 1$ holds
trivially. In
that case, Theorem \ref{lemlowerboundvarselecrow} asserts that
%
%
\begin{eqnarray}
\label{eql0-rowsparse-variable-selection} &&\sup_{\mathcal{P}_{0}
(R_0,\sigma^2)} \E
\bigl\llVert\sin\Theta(\hat{\mathcal{S}}, \mathcal{S}) \bigr\rrVert
_F
\nonumber\\[-8pt]\\[-8pt]
&&\qquad\geq c \biggl\{ (R_0 - d) \frac{\sigma^2}{n} \biggl( 1 + \log
\frac{p-d}{R_q - d} \biggr) \land1 \biggr\}^{1/2}.
\nonumber
\end{eqnarray}
When $q \in(0,2)$ the transition between the $T < \gamma^{q/2}$
and $T \geq\gamma^{q/2}$ regimes involves lower order ($\log
\log$) terms
that can be seen in (\ref{eqvariable-selection-lower-order-terms}).
Under Condition \ref{conrow-sparsity},
(\ref{eqsparse-variable-selection-lower-bound}) can be simplified to
%
%
\begin{eqnarray}
\label{eqlq-rowsparse-variable-selection} && \sup_{\mathcal{P}_{0}
(R_0,\sigma^2)} \E
\bigl\llVert\sin\Theta(\hat{\mathcal{S}}, \mathcal{S}) \bigr\rrVert
_F
\nonumber\\[-8pt]\\[-8pt]
&&\qquad\geq c \biggl\{ (R_q - d) \frac{\sigma^2}{n} \biggl( 1 + (1-a)
\log
\frac{ (p-d)^{1-q/2} }{ R_q - d } \biggr) \land1 \biggr\}^{
{1}/{2} - {q}/{2}}.
\nonumber
\end{eqnarray}

%
%
\begin{theorem}[(Row sparse parameter estimation)]
\label{lemlowerboundparrow}
Let $q \in[0,2)$ and $(p,d,R_q)$ satisfy
\[
2 \leq d \quad\mbox{and}\quad2d \leq R_q - d \leq(p - d)^{1-q/2},
\]
and let $T$ and $\gamma$ be defined as in (\ref{eqgamma-T-definition}).
There exists a universal constant $c > 0$ such that
every estimator $\hat{\mathcal{S}}$ satisfies the following.
If $T < (d \gamma)^{q/2}$, then
%
%
\begin{equation}
\label{eqsparse-estimation-lower-bound} \sup_{\mathcal{P}_{q} (\sigma
^2,R_q)} \E\bigl
\llVert\sin\Theta(\hat{\mathcal{S}}, \mathcal{S}) \bigr\rrVert
_F \geq c \biggl\{ (R_q - d) \biggl(\frac{d \sigma^2}{n}
\biggr)^{1 -{q/2}} \land d \biggr\}^{1/2}.
\end{equation}
Otherwise,
%
%
\begin{equation}
\label{eqdense-estimation-lower-bound} \sup_{\mathcal{P}_{q} (\sigma
^2,R_q)} \E\bigl
\llVert\sin\Theta(\hat{\mathcal{S}}, \mathcal{S}) \bigr\rrVert
_F \geq c \biggl\{ \frac{d(p-d)\sigma^2}{n} \land d \biggr
\}^{1/2}.
\end{equation}
%
\end{theorem}
This result with
(\ref{eql0-rowsparse-variable-selection}) implies Theorem \ref
{thmlowerboundrow0},
and with (\ref{eqlq-rowsparse-variable-selection}) it implies
Theorem~\ref{thmlowerboundrowgeneral}.


%
%
\begin{theorem}[(Column sparse estimation)]
\label{lemlowerboundvarcol}
Let $q \in[0,2)$ and $(p,d,R_q)$ satisfy
\[
4 \leq(p - d)/d \quad\mbox{and}\quad d \leq d (R_q - 1) \leq(p -
d)^{1-q/2},
\]
and recall the definition of $T_*$ in (\ref{eqTstar-definition}).
There exists a universal constant $c > 0$ such that
every estimator $\hat{\mathcal{S}}$ satisfies the following.
If $T_* < \gamma^{q/2}$, then
%
%
\begin{eqnarray}
\label{eqcolumnsparse-lower-bound} && \sup_{\mathcal{P}_{q}^* (\sigma
^2,R_q)} \E\bigl
\llVert\sin\Theta(\hat{\mathcal{S}}, \mathcal{S}) \bigr\rrVert
_F
\nonumber\\[-8pt]\\[-8pt]
&&\qquad\geq c \biggl\{ d (R_q - 1) \biggl[ \frac{\sigma^2}{n} \bigl(
1 -
\log\bigl( T_*/\gamma^{q/2} \bigr) \bigr) \biggr]^{1-
{q/2}} \land d
\biggr\}^{1/2}.
\nonumber
\end{eqnarray}
Otherwise,
%
%
\begin{equation}
\label{eqcolumn-dense-lower-bound} \sup_{\mathcal{P}_{q}^* (\sigma
^2,R_q)} \E\bigl\llVert
\sin\Theta(\hat{\mathcal{S}}, \mathcal{S}) \bigr\rrVert_F \geq
c \biggl\{ \frac{(p-d) \sigma^2}{n} \land d \biggr\}^{1/2}.
\end{equation}
\end{theorem}

In the next section we setup a general technique, using Fano's
inequality and
Stiefel manifold embeddings,
for obtaining minimax lower bounds in principal subspace estimation problems.
Then we move on to proving Theorems \ref
{lemlowerboundvarselecrow} and~\ref{lemlowerboundvarcol}.

\subsection{Lower bounds for principal subspace estimation via Fano method}
\label{subsecfano}
Our main tool for proving minimax lower bounds is the
generalized Fano method. We quote the following version from
\citet{Yu1997}, Lemma 3.

%
%
\begin{lemma}[(Generalized Fano method)]
\label{lemfano}
Let $N \geq1$ be an integer and $\{\theta_1,\ldots,\theta_N \}
\subset\Theta$
index a collection of probability measures $\prob_{\theta_i}$ on a measurable
space $(\mathcal{X}, \mathcal{A})$. Let $d$ be a pseudometric on
$\Theta$ and
suppose that for all $i \neq j$
\[
d(\theta_i, \theta_j) \geq\alpha_N
\]
and, the Kullback--Leibler (KL) divergence
\[
D( {\prob_{\theta_i}} \Vert{\prob_{\theta_j}} ) \leq
\beta_N.
\]
Then every $\mathcal{A}$-measurable estimator $\hat{\theta}$ satisfies
\[
\max_i \E_{\theta_i} d(\hat{\theta},
\theta_i) \geq\frac{\alpha_N}{2} \biggl[ 1 - \frac{\beta_N + \log
2}{\log N}
\biggr].
\]
\end{lemma}

The calculations required for applying Lemma \ref{lemfano} are
tractable when
$\{ \prob_{\theta_i} \}$ is a collection of multivariate\vadjust{\goodbreak} Normal distributions.
Let $A \in\mathbb{V}_{p,d}$ and consider the mean zero $p$-variate Normal
distribution with covariance matrix
%
%
\begin{equation}
\label{eqspiked-covariance-model} \Sigma(A) = b AA^T +
I_{p} = (1+b) A A^T + \bigl(I_{p} - A
A^T\bigr),
\end{equation}
where $b > 0$.
The noise-to-signal ratio of the principal $d$-dimensional subspace of these
covariance matrices is
\[
\sigma^{2} = \frac{1+b}{b^2}
\]
and can choose $b$ to achieve any $\sigma^2 > 0$. The KL divergence between
these multivariate Normal distributions has a simple, exact expression
given in
the following lemma. The proof is straightforward and contained in
Appendix~\ref{secadditional-proof-lower-bound}.

%
%
\begin{lemma}[(KL divergence)]
\label{lemkl-divergence}
For $i = 1,2$, let $A_i \in\mathbb{V}_{p,d}$, $b \geq0$,
\[
\Sigma(A_i) = (1+b) A_i A_i^T
+ \bigl(I_{p} - A_i A_i^T\bigr),
\]
and $\prob_i$ be the $n$-fold product of the $\normal(0,\Sigma(A_i))$
probability measure. Then
\[
D( {\prob_1} \Vert{\prob_2} ) = \frac{n b^2}{1+b}
\bigl\llVert\sin(A_1, A_2) \bigr\rrVert
_F^2.
\]
\end{lemma}

The KL divergence between the probability measures in Lemma \ref
{lemkl-divergence}
is equivalent to the subspace distance. In applying Lemma \ref{lemfano},
we will need to find packing sets in $\mathbb{V}_{p,d}$ that satisfy
the sparsity
constraints of the model \emph{and} have small diameter according to the
subspace Frobenius distance. The next lemma, proved in the \hyperref
[app]{Appendix},
provides a
general method for constructing such local packing sets.

%
%
\begin{lemma}[(Local Stiefel embedding)]
\label{lemstiefel-embedding}
Let $1 \leq k \leq d < p$ and the function $A_\varepsilon\dvtx\mathbb
{V}_{p-d,k}
\mapsto\mathbb{V}_{p,d}$ be defined in block form as
%
%
\begin{equation}
\label{eqstiefel-embedding} A_\varepsilon(J) = \left[\matrix{ \bigl(1-
\varepsilon^2\bigr)^{1/2} I_{k} & 0
\cr
0 &
I_{d-k}
\cr
\varepsilon J & 0}\right]
\end{equation}
for $0 \leq\varepsilon\leq1$.
If $J_1, J_2 \in\mathbb{V}_{p-d,k}$, then
\[
\varepsilon^2 \bigl(1-\varepsilon^2\bigr) \llVert
J_1 - J_2 \rrVert_F^2 \leq
\bigl\llVert\sin\bigl(A_{\varepsilon}(J_1), A_{\varepsilon}(J_2)
\bigr) \bigr\rrVert_F^2 \leq\varepsilon^2
\llVert J_1 - J_2 \rrVert_F^2.
\]
\end{lemma}
This lemma allows us to convert global $O(1)$-separated packing sets in
$\mathbb{V}_{p-d,k}$ into $O(\varepsilon)$-separated packing sets in
$\mathbb{V}_{p,d}$ that are localized within a $O(\varepsilon)$-diameter.
Note that
\[
\llVert J_i - J_j \rrVert_F \leq\llVert
J_i \rrVert_F + \llVert J_j \rrVert
_F \leq2 \sqrt{k}.
\]
By using Lemma \ref{lemstiefel-embedding} in conjunction with Lemmas
\ref{lemfano} and \ref{lemkl-divergence}, we have the following
generic method for
lower bounding the minimax risk of estimating the principal subspace of a
covariance matrix.

%
%
\begin{lemma}
\label{lemgeneric-stiefel-lower-bound}
Let $\varepsilon\in[0, 1]$ and $\{J_1,\ldots,J_N\} \subseteq\mathbb
{V}_{p-d,k}$
for $1 \leq k \leq d < p$.
For each $i = 1,\ldots, N$,
let $\prob_i$ be the $n$-fold product of the $\normal(0, \Sigma
(A_\varepsilon(J_i)))$
probability measure,
where
$\Sigma(\cdot)$ is defined in (\ref{eqspiked-covariance-model})
and
$A_\varepsilon(\cdot)$ is defined in (\ref{eqstiefel-embedding}).
If
\[
\min_{i \neq j} \llVert J_i - J_j
\rrVert_F \geq\delta_N,
\]
then every estimator $\hat{\mathcal{A}}$ of
$\mathcal{A}_i:=\colspan(A_\varepsilon(J_i))$ satisfies
\[
\max_i \E_i \bigl\llVert\sin\Theta(\hat{
\mathcal{A}}, \mathcal{A}_i) \bigr\rrVert_F \geq
\frac{\delta_N \varepsilon\sqrt{1-\varepsilon^2}}{2} \biggl[ 1 - \frac{
4 n k\varepsilon^2 / \sigma^2 + \log2}{\log N} \biggr],
\]
where $\sigma^2 = (1+b) / b^2$.
\end{lemma}

\subsection{Proofs of the main lower bounds}


\mbox{}

\begin{pf*}{Proof of Theorem \ref{lemlowerboundvarselecrow}}
The following lemma, derived from Massart [(\citeyear{Massart2007}),
Lemma 4.10],
allows us to
analyze the variable selection aspect.

%
%
\begin{lemma}[(Hypercube construction)]
\label{lemhypercube}
Let $m$ be an integer satisfying $e \leq m$ and let $s \in[1,m]$.
There exists a subset $\{J_1,\ldots,J_N\} \subseteq\mathbb{V}_{m,1}$
satisfying
the following properties:
\begin{enumerate}
\item
$\llVert J_i \rrVert_{2,0} \leq s$ for all $i$,
\item
$\llVert J_i - J_j \rrVert_2^2 \geq1/4$
for all $i \neq j$, and
\item
$\log N \geq\max\{ c s[ 1 + \log(m/s) ], \log(m) \}$,
where $c > 1/30$ is an absolute constant.
\end{enumerate}
\end{lemma}

%
%
\begin{proposition}\label{prolq-l0-relation}
If $J \in\mathbb{V}_{m,d}$ and $q \in(0,2]$, then
$\llVert J \rrVert_{2,q}^q \leq d^{q/2} \llVert
J \rrVert_{2,0}^{1 - {q/2}}$.
\end{proposition}

Let $\rho\in(0,1]$ and $\{J_1,\ldots, J_N\} \subseteq\mathbb
{V}_{m,1}$ be the
subset given by Lemma \ref{lemhypercube} with $m = p-d$ and $s = \max
\{1, (p-d) \rho\}$.
Then
\begin{eqnarray*}
\log N &\geq& \max\bigl\{ c s\bigl( 1 + \log\bigl[(p-d)/s\bigr] \bigr),
\log(p-d)
\bigr\}
\\
&\geq& \max\bigl\{ (1/30) (p-d)\rho( 1 - \log\rho), \log(p-d) \bigr\}.
\end{eqnarray*}
Applying Lemma \ref{lemgeneric-stiefel-lower-bound}, with $k = 1$,
$\delta_N = 1/2$, and
$b$ chosen so that $(1+b)/b^2 = \sigma^2$, yields
%
%
\begin{eqnarray}\label{eqrowsparse-fano}
&&
\max_i \E_i \bigl\llVert\sin\Theta(\hat{
\mathcal{A}}, \mathcal{A}_i) \bigr\rrVert_F
\nonumber\\
&&\qquad\geq
\frac{\varepsilon}{4 \sqrt{2}} \biggl[ 1 - \frac{4 n \varepsilon^2
/ \sigma^2}{(1/30) (p-d) \rho(1-\log\rho)} - \frac{\log2}{\log(p-d)}
\biggr]
\nonumber\\[-8pt]\\[-8pt]
&&\qquad= \frac{\varepsilon}{4 \sqrt{2}} \biggl[ 1 - \frac{120 \varepsilon
^2}{\gamma\rho(1-\log\rho)} -
\frac{\log2}{\log(p-d)} \biggr]
\nonumber\\
&&\qquad\geq \frac{\varepsilon}{4 \sqrt{2}} \biggl[ \frac{1}{2} -
\frac{120 \varepsilon^2}{\gamma\rho(1-\log\rho)} \biggr]\nonumber
\end{eqnarray}
for every estimator $\hat{\mathcal{A}}$ and all $\varepsilon\in
[0,1/\sqrt{2}]$,
because $p-d \geq4$ by assumption.
Since $J_i \in\mathbb{V}_{p-d,1}$, Proposition \ref
{prolq-l0-relation} implies
%
%
\begin{equation}
\label{eqrow-sparsity-bound} \bigl\llVert A_{\varepsilon}(J_i)
\bigr\rrVert_{2,q} \leq\cases{ d + s, &\quad if $q = 0$\quad and
\cr
\bigl( d +
\varepsilon^q s^{(2-q)/2} \bigr)^{1/q}, &\quad if $0 < q < 2$.}
\end{equation}
For every $q \in[0,2)$
\[
d + \varepsilon^q s^{(2-q)/2} \leq R_q \quad\iff\quad
\varepsilon^{2q} \leq\frac{(R_q - d)^2}{s^{2-q}} = \frac{(R_q -
d)^2}{\max\{1, (p-d)\rho\}^{2-q}}.
\]
Thus, (\ref{eqrow-sparsity-bound}) implies that the constraint
%
%
\begin{equation}
\label{eqrow-sparsity-constraint} \varepsilon^{2q} \leq\min\bigl\{ (T/
\rho)^2 \rho^q, (R_q - d)^2 \bigr
\}
\end{equation}
is sufficient for $\mathcal{A}_i \in\mathcal{M}_{q} (R_q)$
and hence $\prob_i \in\mathcal{P}_{q} (\sigma^2,R_q)$.
Now fix
\[
\varepsilon^2 = \tfrac{1}{480} \gamma\rho(1-\log\rho) \land
\tfrac{1}{2}.
\]
If we can choose $\rho\in(0,1]$ such that
(\ref{eqrow-sparsity-constraint}) is satisfied,
then by (\ref{eqrowsparse-fano}),
\[
\sup_{\mathcal{P}_{q} (\sigma^2,R_q)} \E\bigl\llVert\sin\Theta(\hat
{\mathcal{S}},
\mathcal{S}) \bigr\rrVert_F \geq\max_i
\E_i \bigl\llVert\sin\Theta(\hat{\mathcal{A}},
\mathcal{A}_i) \bigr\rrVert_F \geq\frac{\varepsilon}{16\sqrt{2}}.
\]
%
Choose $\rho\in(0,1]$ to be the unique solution of the equation
%
%
\begin{equation}
\label{eqrho-definition} \rho= \cases{ T \bigl[ \gamma(1- \log\rho)
\bigr]^{-{q}/{2}}, &\quad if $T < \gamma^{q/2}$\quad and
\cr
1, &\quad otherwise.}
\end{equation}
We will verify that $\varepsilon$ and $\rho$ satisfy (\ref
{eqrow-sparsity-constraint}).
The assumption that $1 \leq R_q - d$
guarantees that $\varepsilon^{2q} \leq(R_q - d)^2$, because $\varepsilon
^{2q} \leq1$.
If $T < \gamma^{q/2}$, then
\[
(T/\rho)^2 \rho^q = \bigl[ \gamma\rho(1-\log\rho)
\bigr]^q \geq\varepsilon^{2q}.
\]
If $T \geq\gamma^{q/2}$, then $\rho= 1$ and
\[
(T/\rho)^2 \rho^q = T^2 \geq
\gamma^{q} \geq\varepsilon^{2q}.
\]
Thus, (\ref{eqrow-sparsity-constraint}) holds and so
\[
\sup_{\mathcal{P}_{q} (\sigma^2,R_q)} \E\bigl\llVert\sin\Theta(\hat
{\mathcal{S}},
\mathcal{S}) \bigr\rrVert_F \geq\frac{\varepsilon}{16\sqrt{2}} \geq
\frac{1}{496} \bigl[ \gamma\rho(1 - \log\rho) \bigr]^{1/2} \land
\frac{1}{32}.
\]
%
Now we substitute (\ref{eqrho-definition}) and the definitions of
$\gamma$ and $T$ into the above inequality to get
the following lower bounds.
If $T < \gamma^{q/2}$, then
%
%
\begin{eqnarray}\label{eqvariable-selection-lower-order-terms}
\gamma\rho( 1 - \log\rho)
\nonumber
&=& T \gamma^{1-q/2} \{1 - \log\rho
\}^{1 - {q/2}}
\\
&=& T \gamma^{1-q/2}
\biggl\{1 - \log\bigl(T / \gamma^{q/2} \bigr) + \frac{q}{2}
\log(1-\log\rho) \biggr\}^{1 - {q/2}}
\\
&\geq& T \gamma^{1-q/2} \bigl\{1 - \log\bigl(T /
\gamma^{q/2} \bigr) \bigr\}^{1 - {q/2}}\nonumber
\end{eqnarray}
and so
\begin{eqnarray*}
&& \sup_{\mathcal{P}_{q} (\sigma^2,R_q)} \E\bigl\llVert\sin\Theta(\hat
{\mathcal{S}},
\mathcal{S}) \bigr\rrVert_F
\\
&&\qquad\geq c_0 \biggl\{ (R_q - d) \biggl[ \frac{\sigma^2}{n}
\bigl( 1 - \log\bigl(T / \gamma^{q/2} \bigr) \bigr)
\biggr]^{1-q/2} \land1 \biggr\}^{1/2}.
\end{eqnarray*}
If $T \geq\gamma^{q/2}$,
then $\gamma\rho(1-\log\rho) = \gamma$ and
\[
\sup_{\mathcal{P}_{q} (\sigma^2,R_q)} \E\bigl\llVert\sin\Theta(\hat
{\mathcal{S}},
\mathcal{S}) \bigr\rrVert_F \geq c_0 ( \gamma\land1
)^{1/2} = c_0 \biggl\{ \frac{(p-d) \sigma^2}{n} \land1 \biggr
\}^{1/2}.
\]
\upqed
\end{pf*}

%
\begin{pf*}{Proof of Theorem \ref{lemlowerboundparrow}}
For a fixed subset of $s$ variables, the challenge in estimating the principal
subspace of these variables is captured by the richness of packing sets
in the
Stiefel manifold $\mathbb{V}_{s,d}$. A packing set in the Stiefel
manifold can
be constructed from a packing set in the Grassman manifold by choosing
a single
element of the Stiefel manifold as a representative for each element of the
packing set in the Grassmann manifold. This is well defined, because the
subspace distance is invariant to the choice of basis. The following lemma
specializes known results Pajor [(\citeyear{Pajor1998}), Proposition 8] for
packing sets
in the Grassman manifold.

%
%
\begin{lemma}[{[See \citet{Pajor1998}]}]
\label{lemgrassmann-packing}
Let $k$ and $s$ be integers satisfying $1 \leq k \leq s - k$, and let
$\delta> 0$
There exists a subset $\{J_1,\ldots,J_N\} \subseteq\mathbb{V}_{s,k}$
satisfying
the following properties:
\begin{enumerate}
\item
$\llVert\sin(J_i,J_j) \rrVert_F \geq\sqrt{k} \delta$
for all $i \neq j$, and
\item
$\log N \geq k(s-k) \log(c_2 / \delta)$,
where $c_2 > 0$ is an absolute constant.
\end{enumerate}
\end{lemma}
To apply this result to Lemma \ref{lemgeneric-stiefel-lower-bound} we will
use Proposition \ref{procosine-trace} to
convert the lower bound on the subspace distance into a lower bound on the
Frobenius distance between orthonormal matrices. Thus,
%
%
\begin{equation}
\label{eqcosine-trace-inequality} \llVert J_i - J_j
\rrVert_F \geq\bigl\llVert\sin\Theta(J_i,J_j)
\bigr\rrVert_F \geq\sqrt{k} \delta.
\end{equation}
Let $\rho\in(0,1]$ and $s = \max\{2d, \lfloor(p-d) \rho\rfloor\}$.
Invoke Lemma \ref{lemgrassmann-packing} with $k = d$ and $\delta= c_2/e$,
where $c_2 >0$ is the constant given by Lemma \ref{lemgrassmann-packing}.
Let $\{J_1,\ldots,J_N\} \subseteq\mathbb{V}_{p-d,d}$ be the subset
given by
Lemma \ref{lemgrassmann-packing} after augmenting with rows of zeroes
if necessary.
Then
\[
\log N \geq d(s-d) \geq\max\bigl\{d (s / 2), d^2 \bigr\} \geq\max
\bigl\{(d/4) (p-d) \rho, d^2 \bigr\}
\]
and by (\ref{eqcosine-trace-inequality}),
\[
\llVert J_i - J_j \rrVert_F^2
\geq d (c_2/e)^2
\]
for all $i \neq j$.
The rest of this proof mirrors that of Theorem \ref
{lemlowerboundvarselecrow}.
Let $\varepsilon\in[0,1/\sqrt{2}]$
and apply Lemma \ref{lemgeneric-stiefel-lower-bound} to get
%
%
\begin{eqnarray}\label{eqrowsparse-estimation-fano}
\max_i \E\bigl\llVert\sin\Theta(\hat{\mathcal{A}},
\mathcal{A}_i) \bigr\rrVert_F
&\geq&
\frac{c_2 \sqrt{d} \varepsilon}{2 \sqrt{2} e} \biggl[ 1 - \frac{ 4 n
d\varepsilon^2 / \sigma^2}{(d/4)(p-d)\rho} - \frac{\log2}{d^2} \biggr]
\nonumber\\[-8pt]\\[-8pt]
&\geq& c_1 \sqrt{d} \varepsilon
\biggl[ \frac{1}{2} - \frac{ 16 \varepsilon^2}{\gamma\rho} \biggr
],\nonumber
\end{eqnarray}
where $\gamma$ is defined in (\ref{eqgamma-T-definition})
and we used the assumption that $d \geq2$.
Since $J_i \in\mathbb{V}_{p-d,d}$, Proposition \ref
{prolq-l0-relation} implies
\[
\bigl\llVert A_{\varepsilon}(J_i) \bigr\rrVert_{2,q}
\leq\cases{ d + s, &\quad if $q = 0$\quad and
\cr
\bigl( d + d^{q/2}
\varepsilon^q s^{(2-q)/2} \bigr)^{1/q}, &\quad if $0 < q < 2$.}
\]
For every $q \in(0,2]$
\begin{eqnarray*}
d + d^{q/2} \varepsilon^q s^{(2-q)/2} &\leq&
R_q \quad\iff\\
d^q \varepsilon^{2q} &\leq&
\frac{(R_q - d)^2}{s^{2-q}} = \frac{(R_q - d)^2}{\max\{2d, (p-d)\rho\}^{2-q}}.
\end{eqnarray*}
So $\varepsilon$ and $\rho$ must satisfy the constraint
%
%
\begin{equation}
\label{eqrow-sparsity-constraint-2} d^q \varepsilon^{2q}
\leq\min\biggl\{ (T/\rho)^2 \rho^q, \frac{(R_q - d)^2}{(2d)^{2-q}}
\biggr\}
\end{equation}
to ensure that $\prob_i \in\mathcal{P}_{q} (\sigma^2,R_q)$.
Fix
%
%
\begin{equation}
\label{eqepsilon-definition-2} \varepsilon^2 = \tfrac{1}{64}
\gamma\rho\land\tfrac{1}{2}
\end{equation}
and
%
%
\begin{equation}
\label{eqrho-definition-2} \rho= \cases{ T (d \gamma)^{-{q}/{2}},
&\quad
if $T <
(d \gamma)^{q/2}$\quad and
\cr
1, &\quad otherwise.}
\end{equation}
Since $\varepsilon^2 \leq1/2$,
\[
d^q \varepsilon^{2q} \leq\frac{(R_q - d)^2}{(2d)^{2-q}} \quad\iff\quad
2^q \varepsilon^{2q} \leq\frac{(R_q - d)^2}{4 d^2}
\quad\Longleftarrow\quad
2 d \leq R_q - d,
\]
where the right-hand side is an assumption of the lemma.
That verifies one of the inequalities in (\ref{eqrow-sparsity-constraint-2}).
If $T < (d \gamma)^{q/2}$, then
\[
(T / \rho)^2 \rho^q = (d \gamma\rho)^q
\rho^q \geq d^q \varepsilon^{2q}.
\]
If $T \geq(d \gamma)^{q/2}$, then $\rho= 1$ and
\[
(T / \rho)^2 \rho^q = T^2 \geq(d
\gamma)^q \geq d^q \varepsilon^{2q}.
\]
Thus, (\ref{eqrow-sparsity-constraint-2}) holds and
by (\ref{eqrowsparse-estimation-fano}),
\begin{eqnarray*}
\sup_{\mathcal{P}_{q} (\sigma^2,R_q)} \E\bigl\llVert\sin\Theta(\hat
{\mathcal{S}},
\mathcal{S}) \bigr\rrVert_F &\geq&\max_i
\E_i \bigl\llVert\sin\Theta(\hat{\mathcal{A}},
\mathcal{A}_i) \bigr\rrVert_F
\\
&\geq& c_1 \sqrt{d} \varepsilon\biggl[ \frac{1}{2} -
\frac{ 16 \varepsilon^2}{\gamma^{({2-q})/{q}} \rho} \biggr]
\\
&\geq&\frac{c_1}{4} \sqrt{d} \varepsilon
\\
&\geq& c_0 ( d \gamma\rho\land d )^{1/2}.
\end{eqnarray*}
Finally, we substitute the definition of $\gamma$ and (\ref
{eqrho-definition-2})
into the above inequality to get
the following lower bounds.
If $T < (d \gamma)^{q/2}$, then
\begin{eqnarray*}
\sup_{\mathcal{P}_{q} (\sigma^2,R_q)} \E\bigl\llVert\sin\Theta(\hat
{\mathcal{S}},
\mathcal{S}) \bigr\rrVert_F &\geq& c_0 \bigl\{ T (d
\gamma)^{1 - {q/2}} \land d \bigr\}^{1/2}
\\
&=& c_0 \biggl\{ (R_q - d) \biggl(\frac{d \sigma^2}{n}
\biggr)^{1 - {q/2}} \land d \biggr\}^{1/2}.
\end{eqnarray*}
If $T \geq(d \gamma)^{q/2}$, then
\[
\sup_{\mathcal{P}_{q} (\sigma^2,R_q)} \E\bigl\llVert\sin(\hat{V}, V)
\bigr\rrVert
_F \geq c_0 ( d \gamma\land d )^{1/2} =
c_0 \biggl\{ \frac{d(p-d)\sigma^2}{n} \land d \biggr\}^{1/2}.
\]
\upqed
\end{pf*}

\begin{pf*}{Proof of Theorem \ref{lemlowerboundvarcol}}
The proof is a modification of the proof of Theorem \ref
{lemlowerboundvarselecrow}.
The difficulty of the problem is captured by the difficulty of variable
selection within each column of $V$.
Instead of using a single hypercube construction as in the proof of
Theorem \ref{lemlowerboundvarselecrow}, we apply a hypercube
construction on each
of the $d$ columns. We do this by dividing the $(p-d) \times d$ matrix into
$d$ submatrices of size $\lfloor(p-d)/d \rfloor\times d$, that is,
constructing matrices of the form
\[
\bigl[\matrix{ B_1^T & B_2^T &\cdots& B_d^T & 0 & \cdots}\bigr]^T
\]
and confining the hypercube construction to the $k$th column of each
$\lfloor
(p-d)/d \rfloor\times d$ matrix $B_k$, $k=1,\ldots,d$. This ensures
that the
resulting $(p-d) \times d$ matrix has orthonormal columns with disjoint
supports.

Let $\rho\in(0,1]$ and $s \in\max\{1, \lfloor(p-d)/d \rfloor\rho
\}$.
Applying Lemma \ref{lemhypercube} with $m = \lfloor(p-d)/d \rfloor$,
we obtain a subset $\{J_1,\ldots, J_M\} \subseteq\mathbb{V}_{m,1}$
such that:
\begin{enumerate}
\item$\llVert J_i \rrVert_0 \le s$ for all $i$,
\item$\llVert J_i-J_j \rrVert_2^2 \ge1/4$ for all $i \neq
j$, and
\item$\log M \geq\max\{cs(1+\log(m/s)), \log m \}$,
where $c > 1/30$ is an absolute constant.
\end{enumerate}
Next we will combine the elements of this packing set in $\mathbb
{V}_{m,1}$ to
form a packing set in $\mathbb{V}_{p-d,d}$. A naive approach takes
the $d$-fold product $\{J_1,\ldots,J_M\}^d$, however this results in
too small
a packing distance because two elements of this product set may differ
in only
one column.

We can increase the packing distance by requiring a substantial number of
columns to be different between any two elements of our packing set without
much sacrifice in the size of the final packing set.
This is achieved by applying an additional combinatorial round with the
Gilbert--Varshamov
bound on $M$-ary codes of length $d$ with minimum Hamming distance $d/2$
[\citet{Gilbert52,Varshamov57}].
The $k$th coordinate of each code specifies which element of
$\{J_1,\ldots,J_M\}$ to place in the $k$th column of $B_k$, and so any two
elements of the resulting packing set will differ in at least $d/2$ columns.
Denote the resulting subset of $\mathbb{V}_{p-d,d}$ by~$\mathcal{H}^s$.
We have:
\begin{enumerate}
\item$\llVert H \rrVert_{*,0}\le s$ for all $H\in\mathcal{H}^s$.
\item$\llVert H_1-H_2 \rrVert_2^2\ge d/8$ for all
$H_1,H_2\in\mathcal{H}^s$ such
that $H_1\neq H_2$.
\item$\log N:=\log\llvert\mathcal{H}^s\rrvert
\geq\max\{c d
s(1+\log(m/s)), \log m \}$, where
$c > 0$ is an absolute constant.
\end{enumerate}
Note that the lower bound of $\log m$ in the third item arises by
considering the packing set whose $N$ elements consist of matrices
whose columns in $B_1,\ldots,B_d$ are all equal to some $J_i$ for $i =
1,\ldots, M$. This ensures that $\log N \geq\log M \geq\log m$. From
here, the proof is a straightforward modification of proof of Theorem
\ref{lemlowerboundvarselecrow} with the substitution of $p-d$ by
$(p-d)/d$. For brevity we will only outline the major steps.

Recall the definitions of $T_*$ and $\gamma$ in (\ref{eqTstar-definition}).
Apply Lemma \ref{lemgeneric-stiefel-lower-bound} with the subset
$\mathcal{H}^s$, $k =
d$, $\delta_N = \sqrt{d}/\sqrt{8}$, and $b$ chosen so that
$(1+b)/b^2 =
\sigma^2$. Then
\begin{eqnarray*}
\max_i \E\bigl\llVert\sin\Theta(\hat{\mathcal{A}},
\mathcal{A}_i) \bigr\rrVert_F &\geq& c_0
\sqrt{d} \varepsilon\biggl[ 1 - \frac{4 n \varepsilon^2 / \sigma^2}{c
m \rho(1-\log\rho)} - \frac{\log2}{\log m} \biggr]
\\
&\geq& c_0 \sqrt{d} \varepsilon\biggl[ \frac{1}{4} -
\frac{(8/c) d \varepsilon^2}{\gamma\rho(1-\log\rho)} \biggr]
\end{eqnarray*}
by the assumption that $(p-d)/d \geq4$,
and
\[
\llVert A_i \rrVert
_{*,q} \leq\cases{ 1 + s, &\quad if $q = 0$\quad and
\cr
\bigl( 1 +
\varepsilon^q s^{(2-q)/2} \bigr)^{1/q}, &\quad if $0 < q < 2$.}
\]
The constraint
\[
d^q \varepsilon^{2q} \leq\min\bigl\{ (T_* /
\rho)^2 \rho^q, d^q (R_q -
1)^2 \bigr\}
\]
ensures that $\prob_i \in\mathcal{P}_{q}^* (\sigma^2,R_q)$.
It is satisfied by choosing $\varepsilon$ so that
\[
d \varepsilon^2 = c_1 \gamma\rho(1-\log\rho) \land
\tfrac{1}{2},
\]
where $c_1 > 0$ is a sufficiently small constant,
the assumption that $d <\break  d(R_q - 1)$,
and letting $\rho$ be the unique solution of the equation
\[
\rho= \cases{ T_* \bigl[\gamma(1- \log\rho)\bigr]^{-{q}/{2}},
&\quad if
$T_* <
\gamma^{q/2}$\quad and
\cr
1, &\quad otherwise.}
\]
We conclude that every estimator $\hat{V}$ satisfies
\[
\sup_{\mathcal{P}_{q}^* (\sigma^2,R_q)} \E\bigl\llVert\sin\Theta(\hat
{\mathcal{S}},
\mathcal{S}) \bigr\rrVert_F \geq c_2 \bigl\{ \gamma
\rho(1-\log\rho) \land d \bigr\}^{1/2},
\]
and we have the following explicit lower bounds.
If $T_* < \gamma^{q/2}$, then
\begin{eqnarray*}
&& \sup_{\mathcal{P}_{q}^* (\sigma^2,R_q)} \E\bigl\llVert\sin\Theta(\hat
{\mathcal{S}},
\mathcal{S}) \bigr\rrVert_F
\\
&&\qquad\geq c_3 \biggl\{ d (R_q - 1) \biggl[
\frac{\sigma^2}{n} \bigl( 1 - \log\bigl( T_* / \gamma^{q/2}
\bigr) \bigr)
\biggr]^{1-q/2} \land d \biggr\}^{1/2}.
\end{eqnarray*}
If $T_* \geq\gamma^{q/2}$, then
\[
\sup_{\mathcal{P}_{q}^* (\sigma^2,R_q)} \E\bigl\llVert\sin\Theta(\hat
{\mathcal{S}},
\mathcal{S}) \bigr\rrVert_F \geq c_3 \biggl\{
\frac{(p-d) \sigma^2}{n} \land d \biggr\}^{1/2}.
\]
\upqed
\end{pf*}

\section{Upper bound proofs}
\label{secproofupperbound}

\subsection{Proofs of the main upper bounds}\label
{subsecproofmainupper}\label{secproofupperboundmain}
$\Sigma$ and $\SampleCov$ are both invariant under translations of
$\mu$.
Since our estimators only depend on $X_1,\ldots,X_n$ only through~$\SampleCov$,
we will assume without loss of generality that $\mu= 0$ for the
remainder of
the paper. The sample covariance matrix can be written~as
\[
\SampleCov= \frac{1}{n} \sum_{i=1}^n
(X_i - \bar{X}) (X_i - \bar{X})^T =
\frac{1}{n} \sum_{i=1}^n
X_i X_i^T - \bar{X} \bar{X}^T.
\]
It can be show that $\bar{X} \bar{X}^T$ is a higher order term that is
negligible [see the proofs in \citet{Vu2012}, for an example of such
arguments].
Therefore, we will ignore this term and focus on the dominating
$\frac{1}{n} \sum_{i=1}^n X_i X_i^T$ term in our proofs below.
\begin{pf*}{Proof of Theorem \ref{thmupperboundrowsparse}}
Again, we start from Corollary \ref{lemcurvature-centered}, which gives
\[
\hat{\varepsilon}^2:= \bigl\llVert\sin\Theta(\hat{\mathcal{S}},
\mathcal{S}) \bigr\rrVert_F^2 \leq\frac{\langle\SampleCov-\Sigma,\hat
V\hat V^T-VV^T \rangle
}{\lambda
_d - \lambda_{d+1}}.
\]
To get the correct dependence on $\lambda_i$ and for general
values of $q$, we need a more refined analysis to control
the random variable $\langle\SampleCov-\Sigma,\hat V\hat V^T-VV^T
\rangle$.
Let
\[
W:=\SampleCov- \Sigma,\qquad \Pi:=VV^T\quad\mbox{and}\quad\hat{\Pi}:=
\hat{V} \hat{V}^T.\vadjust{\goodbreak}
\]
Recall that for an orthogonal projector $\Pi$ we write
$\Pi^\perp:=I-\Pi$.
By Proposition~\ref{procurvature-identity} we have
%
%
\begin{eqnarray}
\label{eqdecomposecurvature}
\langle W,\hat\Pi-\Pi\rangle&=&-\bigl\langle W,\Pi\hat\Pi^\bot\Pi\bigr
\rangle+2\bigl\langle W,\Pi^\bot\hat\Pi\Pi\bigr\rangle+\bigl\langle W,
\Pi^\bot\hat\Pi\Pi^\bot\bigr\rangle
\\[-1pt]
\label{equpper-bound-three-terms}
&=:&-T_1+2T_2+T_3.
\end{eqnarray}
We will control $T_1$ (the upper-quadratic term), $T_2$ (the
cross-product term), and $T_3$ (the lower-quadratic term) separately.

\subsubsection*{Controlling $T_1$}
%
%
\begin{eqnarray}
\label{eqcontrolT1a} \llvert
T_1 \rrvert&=& \bigl\llvert\bigl\langle W,\Pi\hat\Pi
^\bot\Pi\bigr\rangle\bigr\rrvert= \bigl\llvert\bigl\langle\Pi W\Pi,
\Pi\hat\Pi^\bot\Pi\bigr\rangle\bigr\rrvert
\nonumber\\[-1pt]
&\le& \llVert\Pi W\Pi\rrVert_2\bigl\llVert\Pi\hat\Pi
^\bot\Pi\bigr\rrVert_{*}= \llVert\Pi W\Pi\rrVert
_2\bigl\llVert\Pi\hat\Pi^\bot\hat\Pi^\bot\Pi
\bigr\rrVert_{*}
\\[-1pt]
&=&\llVert\Pi W\Pi\rrVert_2\bigl\llVert\Pi\hat\Pi^\bot
\bigr\rrVert_{F}^2 \le\llVert\Pi W\Pi\rrVert
_2 \hat\varepsilon^2,
\nonumber
\end{eqnarray}
where $\llVert\cdot\rrVert_{*}$ is the nuclear norm
($\ell_1$ norm
of the singular values) and $\llVert\cdot\rrVert_2$ is
the spectral norm (or operator
norm). By Lemma \ref{lemupper-quadratic}, we have (recall that
we assume $\llVert Z \rrVert_{\psi_2} \leq1$ and $\varepsilon
_n\le1$
for simplicity)
%
%
\begin{equation}
\bigl\llVert{\llVert\Pi W\Pi\rrVert_2}\bigr\rrVert
_{\psi
_1}\le c_1 \lambda_1 \sqrt{d/n},
\end{equation}
where $c_1$ is a universal constant.
Define
\[
\Omega_1= \biggl\{\llvert T_1 \rrvert\ge
c_1\sqrt{\frac{d}{n}} \log n \lambda_1\hat
\varepsilon^2 \biggr\}.
\]
Then, when $n\ge2$ we have
%
%
\begin{equation}
\label{eqOmega1bound} \Pr(\Omega_1)\le
\Pr\bigl(\llVert\Pi W\Pi\rrVert_2\ge c_1
\lambda_1\log n \sqrt{d/n} \bigr)\le(n-1 )^{-1}.
\end{equation}

\subsubsection*{Controlling $T_2$}
%
%
\begin{eqnarray}\label{eqboundT2a}
T_2&=&\bigl\langle W,\Pi^\bot\hat\Pi\Pi\bigr\rangle=\bigl
\langle\Pi^\bot W \Pi,\Pi^\bot\hat\Pi\bigr\rangle
\nonumber\\[-9pt]\\[-9pt]
&\le&\bigl\llVert\Pi^\bot W\Pi\bigr\rrVert_{2,\infty} \bigl
\llVert\Pi^\bot\hat\Pi\bigr\rrVert_{2,1}.
\nonumber
\end{eqnarray}
To bound $\llVert\Pi^\bot\hat\Pi\rrVert_{2,1}$,
let the rows of $\Pi^\perp\hat\Pi$ be denoted
by $\phi_1,\ldots, \phi_p$ and $t > 0$.
Using a standard argument of bounding $\ell_1$ norm
by the $\ell_q$ and $\ell_2$ norms [e.g., \citet{Raskutti2011},
Lemma 5], we have for all $t>0$, $0<q\le1$,
%
%
\begin{eqnarray}
\label{eq2qboundproductterm}
\bigl\llVert\Pi^\perp\hat\Pi\bigr\rrVert_{2,1} &=& \sum
_{i=1}^p \llVert\phi_i
\rrVert_2
\nonumber\\[-1pt]
&\leq&\Biggl[ \sum_{i=1}^p \llVert
\phi_i \rrVert_2^q \Biggr]^{1/2}
\Biggl[ \sum_{i=1}^p \llVert
\phi_i \rrVert_2^2 \Biggr]^{1/2}
t^{-q/2} + \Biggl[\sum_{i=1}^p
\llVert\phi_i \rrVert_2^q \Biggr]
t^{1-q}
\nonumber\\[-9pt]\\[-9pt]
&=& \bigl\llVert\Pi^\perp\hat\Pi\bigr\rrVert_{2,q}^{q/2}
\bigl\llVert\Pi^\perp\hat\Pi\bigr\rrVert_F
t^{-q/2} + \bigl\llVert\Pi^\perp\hat\Pi\bigr\rrVert
_{2,q}^q t^{1-q}
\nonumber
\\
&\leq&\sqrt{2} R_q^{1/2} t^{-q/2}\hat\varepsilon+ 2
R_q t^{1-q},
\nonumber
\end{eqnarray}
where the last step uses the fact that
\begin{eqnarray*}
\bigl\llVert\Pi^\perp\hat\Pi\bigr\rrVert_{2,q}^q&=&
\bigl\llVert\Pi^\perp\hat V \bigr\rrVert_{2,q}^q
= \llVert\hat V - \Pi\hat V \rrVert_{2,q}^q \leq\llVert
\hat V \rrVert_{2,q}^q + \bigl\llVert V V^T \hat
V \bigr\rrVert_{2,q}^q
\\
&\leq&\llVert\hat V \rrVert_{2,q}^q + \llVert V \rrVert
_{2,q}^q \leq2 R_q.
\end{eqnarray*}
Combining (\ref{eqboundT2a}) and (\ref{eq2qboundproductterm})
we obtain, for all $t>0$, $0<q<1$,
%
%
\begin{equation}
\label{eqboundT2b} T_2\le
\bigl\llVert\Pi^\bot W\Pi\bigr\rrVert_{2,\infty} \bigl(\sqrt{2}
R_q^{1/2} t^{-q/2}\hat\varepsilon+ 2 R_q
t^{1-q} \bigr).
\end{equation}
The case where $q=0$ is simpler and omitted.
Now define
\begin{eqnarray*}
\Omega_2&:=& \bigl\{ T_2\ge20 \bigl( \sqrt{
\lambda_1\lambda_{d+1}}^{1-q/2}(\lambda_d-
\lambda_{d+1})^{q/2} \varepsilon_n\hat\varepsilon
\\
&&\hspace*{24.5pt}\hspace*{16.7pt}{} +\sqrt{\lambda_1\lambda_{d+1}}^{2-q} (
\lambda_{d}-\lambda_{d+1})^{-(1-q)}
\varepsilon_n^2 \bigr) \bigr\}
\\
&=& \bigl\{ T_2\ge t_{2,1} \bigl(\sqrt{2R_q}t_{2,2}^{-q/2}
\hat\varepsilon+2R_qt_{2,2}^{1-q} \bigr) \bigr\},
\\
t_{2,1} &=& 20 \sqrt{\lambda_1\lambda_{d+1}}
\sqrt{\frac{d+\log p}{n}},
\\
t_{2,2} &=& \frac{\sqrt{\lambda_1\lambda_{d+1}}}{\lambda_d-\lambda
_{d+1}} \sqrt{\frac{d+\log p}{n}}.
\end{eqnarray*}
Taking $t=t_{2,2}$ in (\ref{eqboundT2b}) and using
the tail bound result in Lemma \ref{lemcrosstermtail}, we have
%
%
\begin{eqnarray}
\label{eqOmega2bound} \Pr(\Omega_2)&\le&
\Pr\bigl(\bigl\llVert\Pi^{\perp}W\Pi\bigr\rrVert_{2,\infty} \ge
t_{2,1}\bigr)
\nonumber\\
&\le& 2 p5^d\exp\biggl( - \frac{t_{2,1}^2/8}{2 \lambda_1\lambda_{d+1}/n
+ t_{2,1}
\sqrt{\lambda_1\lambda_{d+1}} / n} \biggr)
\\
&\le& p^{-1}.
\nonumber
\end{eqnarray}

\subsubsection*{Controlling $T_3$}

The bound on $T_3$ involves a quadratic form empirical process
over a random set.
Let $\varepsilon\geq0$ and define
\[
\phi(R_q, \varepsilon):= \sup\bigl\{ \bigl\langle W,\Pi^\perp
U U^T \Pi^\perp\bigr\rangle\dvtx U \in\mathbb{V}_{p,d},
\llVert U \rrVert_{2,q}^q \leq R_q, \bigl
\llVert\Pi^\perp U \bigr\rrVert_F \leq\varepsilon\bigr\}.
\]
Then by Lemma \ref{lemlower-quadratic}, we have, with some universal constants
$c_{3}$,
for $x>0$
\[
\Pr\bigl(\phi(R_q,\varepsilon)\ge c_{3}x
\lambda_{d+1} \bigl(\varepsilon_n\varepsilon^2 +
\varepsilon_n^2\varepsilon+\varepsilon_n^4
\bigr) \bigr) \le2\exp\bigl(-x^{2/5}\bigr).
\]
Let $T_3(U)=\langle W,\Pi^\perp UU^T\Pi^\perp\rangle$,
for all $U\in\mathcal{U}_{p} (R_q)$,
where
\[
\mathcal{U}_{p} (R_q):= \bigl\{ U \in
\mathbb{V}_{p,d}\dvtx\colspan(U) \in\mathcal{M}_{p}
(R_q) \bigr\}.
\]
Define function
$g(\varepsilon)=\varepsilon_n\varepsilon^2+\varepsilon_n^2\varepsilon
+\varepsilon_n^4$.
Then for all $\varepsilon\ge0$, we have
$g(\varepsilon)\ge\varepsilon_n^4\ge4d^3/n^2$. On the other hand,
if $\varepsilon=\llVert{\sin\Theta}(U,V) \rrVert_F$, then
$\varepsilon^2\le2d$ and hence\vspace*{1pt} $g(\varepsilon)\le
g(\sqrt{2d}) =2d+\sqrt{2d}+1$.
Let $\mu=\varepsilon_n^4$ and $J=\lceil\log_2 (g(\sqrt{2d})/\mu
)\rceil$.
Then we have $J\le3\log n+6/5$.

Note that $g$ is strictly increasing
on $[0,\sqrt{2d}]$. Then
we have the following peeling argument:
\begin{eqnarray*}
&&\Pr\bigl[\exists U\in\mathcal{U}_{p} (R_q)\dvtx
T_3(U)\ge2c_3\lambda_{d+1}(\log
n)^{5/2} g\bigl(\bigl\llVert\sin(U,V) \bigr\rrVert_F
\bigr) \bigr]
\\
&&\qquad\le\Pr\bigl[ \exists1\le j\le J, U\in\mathcal{U}_{p}
(R_q)\dvtx2^{j-1}\mu\le g\bigl(\bigl\llVert\sin\Theta(U,V)
\bigr\rrVert_F\bigr)\le2^j\mu,
\\
&&\hspace*{115.2pt} T_3(U)\ge2c_3\lambda_{d+1}(\log
n)^{5/2}g\bigl(\bigl\llVert\sin\Theta(U,V) \bigr\rrVert
_F\bigr) \bigr]
\\
&&\qquad\le\sum_{j=1}^J \Pr\bigl[\phi
\bigl(R_q,g^{-1}\bigl(2^j\mu\bigr)\bigr)\ge
c_3\lambda_{d+1}(\log n)^{5/2}2^{j}\mu
\bigr]
\\
&&\qquad\le J 2 n^{-1}\le\frac{6\log n}{n}+\frac{3}{n}.
\end{eqnarray*}
Define
\[
\Omega_3:= \bigl\{ \phi(R_q,\hat\varepsilon) \ge
c_3(\log n)^{5/2}\lambda_{d+1}\bigl(
\varepsilon_n\hat\varepsilon^2+\varepsilon_n^2
\hat\varepsilon+\varepsilon_n^4\bigr) \bigr\}.
\]
Then we have proved that
\[
\Pr(\Omega_3)\le\frac{6\log n}{n}+\frac{3}{n}.
\]

\subsubsection*{Putting things together}

Now recall the conditions in
(\ref{eqsmallepsilon1}) to (\ref{eqsmallepsilon4}). On $\Omega
_1^c\cap\Omega_2^c \cap\Omega_3^c$,
we have, from (\ref{eqdecomposecurvature}) that
\begin{eqnarray*}
\quad\qquad(\lambda_d-\lambda_{d+1})\hat\varepsilon^2 &\le&
\biggl(c_1\sqrt{\frac{d}{n}}\log n \lambda_1+
c_3\varepsilon_n(\log n)^{5/2}
\lambda_{d+1} \biggr)\hat\varepsilon^2
\\
&&{} +41 \sqrt{\lambda_1\lambda_{d+1}}^{1-q/2}(
\lambda_d-\lambda_{d+1})^{q/2}
\varepsilon_n\hat\varepsilon
\\
&&{} +41 \sqrt{\lambda_1\lambda_{d+1}}^{2-q} (
\lambda_{d}-\lambda_{d+1})^{-(1-q)}
\varepsilon_n^2\quad\Longrightarrow
\\
\frac{1}{2}(\lambda_d-\lambda_{d-1})
\hat\varepsilon^2 &\le& 41 \sqrt{\lambda_1
\lambda_{d+1}}^{1-q/2}(\lambda_d-
\lambda_{d+1})^{q/2} \varepsilon_n\hat\varepsilon
\\
&&{} +41 \sqrt{\lambda_1\lambda_{d+1}}^{2-q} (
\lambda_{d}-\lambda_{d+1})^{-(1-q)}
\varepsilon_n^2\quad\Longrightarrow
\\
\hat\varepsilon&\le&83 \biggl(\frac{\sqrt{\lambda
_1\lambda_{d+1}}} {
\lambda_d-\lambda_{d+1}}
\biggr)^{1-q/2}\varepsilon_n.\hspace*{129pt}\qed
\end{eqnarray*}
\noqed
\end{pf*}

\section{Additional proofs}
\label{secauxproofs}


\begin{pf*}{Proof of Proposition \ref{procosine-trace}}
Let $\gamma_i$ be the cosine of the $i$th canonical angle between the
subspaces spanned by
$V_1$ and $V_2$. By Theorem II.4.11 of \citet{StewartAndSun},
\[
\inf_{Q \in\mathbb{V}_{k,k}} \llVert V_1 - V_2 Q
\rrVert_F^2 = 2 \sum_i
(1 - \gamma_i).
\]
The inequalities
\[
1-x \leq\bigl(1-x^2\bigr) \leq2(1-x)
\]
hold for all $x \in[0,1]$. So
\[
\frac{1}{2} \inf_{Q \in\mathbb{V}_{k,k}} \llVert V_1 -
V_2 Q \rrVert_F^2 \leq\sum
_i \bigl(1-\gamma_i^2\bigr) \leq
\inf_{Q \in\mathbb{V}_{k,k}} \llVert V_1 - V_2 Q \rrVert
_F^2.
\]
Apply the trigonometric identity $\sin^2\theta= 1 -
\cos^2\theta$ to the preceding display to conclude the proof.
\end{pf*}

\subsection{Proofs related to the lower bounds}\label
{secadditional-proof-lower-bound}

\mbox{}

\begin{pf*}{Proof of Lemma \ref{lemkl-divergence}}
Write $\Sigma_i = \Sigma(A_i)$ for $i = 1,2$.
Since $\Sigma_1$ and $\Sigma_2$ are nonsingular and have the same
determinant,
\begin{eqnarray*}
D( {\prob_1} \Vert{\prob_2} ) &=& n D\bigl( {\normal(0,
\Sigma_1)} \Vert{\normal(0, \Sigma_2)} \bigr)
\\
&=& \frac{n}{2} \bigl\{ \tr\bigl(\Sigma_2^{-1}
\Sigma_1\bigr) - p - \log\det\bigl(\Sigma_2^{-1}
\Sigma_1\bigr) \bigr\}
\\
&=& \frac{n}{2} \tr\bigl(\Sigma_2^{-1} (
\Sigma_1 - \Sigma_2) \bigr).
\end{eqnarray*}
Now
\[
\Sigma_2^{-1} = (1 + b)^{-1} A_2
A_2^T + \bigl(I_{p} - A_2
A_2^T\bigr)
\]
and
\[
\Sigma_1 - \Sigma_2 = b \bigl( A_1
A_1^T - A_2 A_2^T
\bigr).
\]
Thus,
\begin{eqnarray*}
&& \tr\bigl(\Sigma_2^{-1} (\Sigma_1 -
\Sigma_2) \bigr)
\\
&&\qquad= \frac{b}{1+b} \bigl\{ (1+b) \bigl\langle I_{p} -
A_2 A_2^T,A_1
A_1^T \bigr\rangle- \bigl\langle A_2
A_2^T,A_2 A_2^T -
A_1 A_1^T \bigr\rangle\bigr\}
\\
&&\qquad= \frac{b-1}{b} \bigl\{ b \bigl\langle I_{p} - A_2
A_2^T,A_1 A_1^T
\bigr\rangle- \bigl\langle I_{p},A_2 A_2^T
- A_2 A_2^T A_1
A_1^T \bigr\rangle\bigr\}
\\
&&\qquad= \frac{b}{1+b} \bigl\{ (1+b) \bigl\langle I_{p} -
A_2 A_2^T,A_1
A_1^T \bigr\rangle- \bigl\langle A_2
A_2^T,I_{p} - A_1
A_1^T \bigr\rangle\bigr\}
\\
&&\qquad= \frac{b^2}{1+b} \bigl\llVert\sin(A_1, A_2) \bigr
\rrVert_F^2
\end{eqnarray*}
by Proposition \ref{procanonical-angles}.
\end{pf*}

%
\begin{pf*}{Proof of Lemma \ref{lemstiefel-embedding}}
By Proposition \ref{procanonical-angles} and the definition of
$A_\varepsilon(\cdot)$,
\begin{eqnarray*}
\bigl\llVert\sin\bigl(A_{\varepsilon}(J_1), A_{\varepsilon}(J_2)
\bigr) \bigr\rrVert_F^2 &=& \frac{1}{2} \bigl
\llVert\bigl[A_{\varepsilon}(J_1)\bigr] \bigl[A_{\varepsilon}(J_1)
\bigr]^T - \bigl[A_{\varepsilon}(J_2)\bigr]
\bigl[A_{\varepsilon}(J_2)\bigr]^T \bigr\rrVert
_F^2
\\
&=& \varepsilon^2\bigl(1-\varepsilon^2\bigr) \llVert
J_1 - J_2 \rrVert_F^2 +
\frac{\varepsilon^4}{2} \bigl\llVert J_1 J_1^T -
J_2 J_2^T \bigr\rrVert_F^2
\\
&\geq&\varepsilon^2 \bigl(1-\varepsilon^2\bigr) \llVert
J_1 - J_2 \rrVert_F^2.
\end{eqnarray*}
The upper bound follows from Proposition \ref{procosine-trace}:
\[
\bigl\llVert\sin\bigl(A_{\varepsilon}(J_1), A_{\varepsilon}(J_2)
\bigr) \bigr\rrVert_F^2 \leq\bigl\llVert
A_{\varepsilon}(J_1) - A_{\varepsilon}(J_2) \bigr
\rrVert_F^2 = \varepsilon^2 \llVert
J_1 - J_2 \rrVert_F^2.
\]
\upqed
\end{pf*}

%
\begin{pf*}{Proof of Lemma \ref{lemhypercube}}
Let $s_0 = \lfloor\min(m / e, s ) \rfloor$.
The assumptions that $m/e \geq1$ and $s \geq1$ guarantee that $s_0
\geq1$.
According to Massart [(\citeyear{Massart2007}), Lem\-ma~4.10] [with
$\alpha= 7/8$ and
$\beta= 8/ (7 e)$],
there exists a subset
$\Omega_m^{s_0} \subseteq\{0,1\}^m$ satisfying the following properties:
\begin{enumerate}
\item
$\llVert\omega\rrVert_0 = s_0$ for all $\omega\in\Omega
_m^{s_0}$,
\item
$\llVert\omega- \omega' \rrVert_0 > s_0 / 4$
for all distinct pairs $\omega, \omega' \in\Omega_m^{s_0}$, and
\item
$\log\llvert\Omega_m^{s_0}\rrvert\geq c s_0 \log
(m/s_0)$, where $c > 0.251$.
\end{enumerate}
Let
\[
\{J_1,\ldots, J_N\}:= \bigl\{ s_0^{-1/2}
\omega\dvtx\omega\in\Omega_m^{s_0} \bigr\}.
\]
Clearly, $\{J_1,\ldots, J_N\} \subseteq\mathbb{V}_{m,1}$ and
\[
\llVert J_i \rrVert_{(2,0)} = \llVert\omega\rrVert
_0 = s_0 \leq s
\]
for every $i$.
If $i \neq j$, then
\[
\llVert J_i - J_j \rrVert_F^2
= s_0^{-1} \llVert\omega_i -
\omega_j \rrVert_0 > 1/4.
\]
The cardinality of $\{J_1,\ldots, J_N\}$ satisfies
\[
\log N = \log\bigl\llvert\Omega_m^{s_0}\bigr\rrvert\geq
c s_0 \log(m / s_0).
\]
As a function of $s_0$, the above right-hand side is increasing on
the interval $[0, m/e]$. Since $\min(m/e, s) / 2 \leq s_0$ belongs to
that interval
\begin{eqnarray*}
\log N &\geq& c \bigl(\min(m/e, s)/2\bigr) \log\bigl[ m / \bigl(\min
(m/e, s)/2
\bigr) \bigr]
\\
&\geq&(c/2) \min(m/e, s) \log\bigl[ m / \min(m/e, s) \bigr].
\end{eqnarray*}
It is easy to see that
\[
\min(m/e, s) \log\bigl[ m / \min(m/e, s) \bigr] \geq\max\bigl\{ s \log(m/s),
s/e \bigr\}
\]
for all $s \in[1,m]$.
Thus,
\[
\min(m/e, s) \log[ m / \min(m/e, s) \geq(1+e)^{-1} s +
(1+e)^{-1} s \log(m/s)
\]
and
%
%
\begin{equation}
\label{eqnegative-binomial-bound} \log N \geq(c/2) (1+e)^{-1} s
\bigl(1 + \log(m/s)\bigr),
\end{equation}
where $(c/2)(1+e)^{-1} > 1/30$.
If the above right-hand side is $\leq\log m$, then we may repeat the entire
argument from the beginning with $\{J_1,\ldots,J_N\}$ taken to be the
$N = m$ vectors
$\{(1,0,\ldots,0), (0,1,0,\ldots,0),\ldots, (0,\ldots, 0, 1) \}
\subseteq\{0,1\}^m$.
That yields, in combination with
(\ref{eqnegative-binomial-bound}),
\[
\log N \geq\max\bigl\{ (1/30) s\bigl[1+\log(m/s)\bigr], \log m \bigr\}.
\]
\upqed
\end{pf*}

\subsection{Proofs related to the upper bounds}\label
{subsecadditionalproofupper}

\mbox{}

\begin{pf*}{Proof of Lemma \ref{lemcurvature-lemma}}
For brevity, denote the eigenvalues of $A$ by $\lambda_d:=
\lambda_d(A)$.
Let $A = \sum_{i=1}^p \lambda_i u_i u_i^T$ be the spectral
decomposition of $A$
so that $E = \sum_{i=1}^d u_i u_i^T$ and $E^\perp= \sum_{i=d+1}^p
u_i u_i^T$.
Then
\begin{eqnarray*}
\langle A,E - F \rangle&=& \bigl\langle A,E (I_{} - F) -
(I_{} - E) F \bigr\rangle
\\
&=& \bigl\langle E A,F^\perp\bigr\rangle- \bigl\langle E^\perp
A,F \bigr\rangle
\\
&=& \sum_{i=1}^d \lambda_i
\bigl\langle u_i u_i^T,F^\perp\bigr
\rangle- \sum_{i=d+1}^p
\lambda_i \bigl\langle u_i u_i^T,F
\bigr\rangle
\\
&\geq&\lambda_d \sum_{i=1}^d
\bigl\langle u_i u_i^T,F^\perp
\bigr\rangle- \lambda_{d+1} \sum_{i=d+1}^p
\bigl\langle u_i u_i^T,F \bigr\rangle
\\
&=& \lambda_d \bigl\langle E,F^\perp\bigr\rangle-
\lambda_{d+1} \bigl\langle E^\perp,F \bigr\rangle.
\end{eqnarray*}
Since orthogonal projectors are idempotent,
\begin{eqnarray*}
\lambda_d \bigl\langle E,F^\perp\bigr\rangle-
\lambda_{d+1} \bigl\langle E^\perp,F \bigr\rangle&=&
\lambda_d \bigl\langle E F^\perp,E F^\perp\bigr
\rangle- \lambda_{d+1} \bigl\langle E^\perp F,E^\perp F
\bigr\rangle
\\
&=& \lambda_d \bigl\llVert E F^\perp\bigr\rrVert
_F^2 - \lambda_{d+1} \bigl\llVert
E^\perp F \bigr\rrVert_F^2.
\end{eqnarray*}
Now apply Proposition \ref{procanonical-angles} to conclude that
\[
\lambda_d \bigl\llVert E F^\perp\bigr\rrVert
_F^2 - \lambda_{d+1} \bigl\llVert
E^\perp F \bigr\rrVert_F^2 = (
\lambda_d - \lambda_{d+1}) \bigl\llVert\sin\Theta(
\mathcal{E}, \mathcal{F}) \bigr\rrVert_F^2.
\]
\upqed
\end{pf*}

%
%
\begin{proposition}
\label{procurvature-identity}
If $W$ is symmetric, and $E$ and $F$ are orthogonal projectors,
then
%
%
\begin{equation}
\langle W,F - E \rangle= \bigl\langle E^\perp W E^\perp,F
\bigr\rangle- \bigl\langle E W E,F^\perp\bigr\rangle+ 2 \bigl\langle
E^\perp W E,F \bigr\rangle.
\end{equation}
\end{proposition}
\begin{pf}
Using the expansion
\[
W = E^\perp W E^\perp+ E W E + E W E^\perp+
E^\perp W E
\]
and the symmetry of $W$, $F$ and $E$,
we can write
\begin{eqnarray*}
\langle W,F - E \rangle&=& \bigl\langle E^\perp W E^\perp,F -
E \bigr\rangle+ \langle E W E,F - E \rangle
\\
&&{} + 2 \bigl\langle E^\perp W E,F - E \bigr\rangle
\\
&=& \bigl\langle E^\perp W E^\perp,E^\perp(F - E)
\bigr\rangle+ \bigl\langle E W E,E(F - E) \bigr\rangle
\\
&&{} + 2 \bigl\langle E^\perp W E,E^\perp(F - E) \bigr\rangle
\\
&=& \bigl\langle E^\perp W E^\perp,F \bigr\rangle+ \bigl
\langle E W E,E(F - E) \bigr\rangle+ 2 \bigl\langle E^\perp W E,F \bigr
\rangle.
\end{eqnarray*}
Now note that
\[
E(F - E) = EF - E = - E F^\perp.
\]
\upqed
\end{pf}

\section{Empirical process related proofs}
\label{secempiricalprocessproofs}
\subsection{The cross-product term}
This section is dedicated to proving the following bound on the
cross-product term.

%
%
\begin{lemma}
\label{lemcrosstermtail}
There exists a universal constant $c > 0$ such that
\[
\Pr\bigl(\bigl\llVert\Pi^{\perp}W\Pi\bigr\rrVert_{2,\infty}> t
\bigr) \le2p5^d \exp\biggl( - \frac{t^2/8}{2 \lambda_1\lambda_{d+1}/n
+ t
\sqrt{\lambda_1\lambda_{d+1}} / n} \biggr).
\]
\end{lemma}

The proof of Lemma \ref{lemcrosstermtail} builds on the following
two lemmas.
They are adapted from Lemmas 2.2.10 and 2.2.11 of \citet
{vanderVaartAndWellner}.
%
%
\begin{lemma}[(Bernstein's inequality)] \label{lembernstein}
Let $Y_1,\ldots,Y_n$ be independent random variables with zero mean.
Then
\[
\prob\Biggl(\Biggl\llvert{\sum_{i=1}^n
Y_i}\Biggr\rrvert> t \Biggr) \leq2 \exp\biggl( -
\frac{t^2 / 2}{2 \sum_{i=1}^n \llVert Y_i \rrVert_{\psi
_1}^2 + t
\max_{i \leq n} \llVert Y_i \rrVert_{\psi_1}} \biggr)
\]
\end{lemma}

%
%
\begin{lemma}[(Maximal inequality)] \label{lembernstein-maximal}
Let $Y_1,\ldots,Y_m$ be arbitrary random variables that satisfy the bound
\[
\prob\bigl(\llvert Y_i \rrvert> t \bigr) \leq2 \exp\biggl( -
\frac{t^2 / 2}{b + a t} \biggr)
\]
for all $t > 0$ (and $i$) and fixed $a,b > 0$. Then
\[
\Bigl\llVert{ \max_{1 \leq i \leq m} Y_i}\Bigr\rrVert
_{\psi_1} \leq c \bigl( a \log(1+m) + \sqrt{b \log(1+m)} \bigr)
\]
for a universal constant $c > 0$.
\end{lemma}

We bound $\llVert\Pi^\perp(\SampleCov- \Sigma) \Pi\rrVert_{2,\infty}$ by
a standard $\delta$-net argument.
%
%
\begin{proposition}
\label{procovering-argument}
Let $A$ be a $p \times d$ matrix,
$(e_1,\ldots,e_p)$ be the canonical basis of $\Real^p$
and $\mathcal{N}_{\delta}$ be a $\delta$-net of $\Sphere_2^{d-1}$
for some
$\delta\in[0,1)$.
Then
\[
\llVert A \rrVert_{2,\infty} \leq(1-\delta)^{-1} \max
_{1 \leq j \leq p} \max_{u \in\mathcal{N}_{\delta}} \langle e_j,A
u \rangle.
\]
\end{proposition}
\begin{pf}
By duality and compactness, there exists $u_* \in\Sphere^{d-1}$ and
$u \in\mathcal{N}_\delta$ such that
\[
\llVert A \rrVert_{2,\infty} = \max_{1\leq j\leq p} \bigl\llVert
e_j^T A \bigr\rrVert_2
= \max_{1\leq j\leq p} \langle e_j,A u_* \rangle
\]
and $\llVert u_* - u \rrVert_2 \leq\delta$.
Then by the Cauchy--Schwarz inequality,
\begin{eqnarray*}
\llVert A \rrVert_{2,\infty} &=& \max_{1\leq j\leq p} \langle
e_j,A u \rangle+ \bigl\langle e_j,A(u_* - u) \bigr
\rangle
\\
&\leq&\max_{1\leq j\leq p} \langle e_j,A u \rangle+ \delta
\bigl\llVert e_j^T A \bigr\rrVert_2
\\
&\leq&\max_{1\leq j\leq p} \max_{u \in\mathcal{N}_\delta} \langle
e_j,A u \rangle+ \delta\llVert A \rrVert_{2,\infty}.
\end{eqnarray*}
Thus,
\[
\llVert A \rrVert_{2,\infty} \leq(1-\delta)^{-1} \max
_{1 \leq j \leq p} \max_{u \in\mathcal{N}_{\delta}} \langle e_j,A
u \rangle.
\]
\upqed
\end{pf}

The following bound on the covering number of the sphere is well known
[see, e.g., \citet{Ledoux}, Lemma 3.18].
%
%
\begin{proposition}\label{prosphere-covering}
Let $\mathcal{N}_\delta$ be a minimal $\delta$-net of $\Sphere_2^{d-1}$
for $\delta\in(0,1)$. Then
\[
\llvert\mathcal{N}_\delta\rrvert\leq( 1 + 2/\delta)^d.
\]
\end{proposition}

%
%
\begin{proposition}
\label{proorlicz-cauchy-schwarz}
Let $X$ and $Y$ be random variables. Then
\[
\llVert X Y \rrVert_{\psi_1} \leq\llVert X \rrVert_{\psi_2}
\llVert Y \rrVert_{\psi_2}.
\]
\end{proposition}
\begin{pf}
Let $A = X / \llVert X \rrVert_{\psi_2}$ and $Y / \llVert Y \rrVert
_{\psi_2}$.
Using the elementary inequality
\[
\llvert a b \rrvert\leq\tfrac{1}{2}\bigl(a^2 + b^2
\bigr)
\]
and the triangle inequality we have that
\[
\llVert A B \rrVert_{\psi_1} \leq\tfrac{1}{2} \bigl( \bigl\llVert
A^2 \bigr\rrVert_{\psi_1} + \bigl\llVert B^2
\bigr\rrVert_{\psi_1} \bigr) = \tfrac{1}{2} \bigl( \llVert A \rrVert
_{\psi_2}^2 + \llVert B \rrVert_{\psi_2}^2
\bigr) =1.
\]
Multiplying both sides of the inequality by $\llVert X \rrVert_{\psi
_2} \llVert Y \rrVert_{\psi_2}$
gives the desired result.
\end{pf}

\begin{pf*}{Proof of Lemma \ref{lemcrosstermtail}}
Let $N_\delta$ be a
minimal $\delta$-net in $\Sphere_2^{d-1}$
for some $\delta\in(0,1)$ to be chosen later.
By Proposition \ref{procovering-argument} we have
\[
\bigl\llVert\Pi^\bot W\Pi\bigr\rrVert_{2,\infty} \le
\frac{1}{1-\delta}\max_{1\le j\le p}\max_{u\in N_\delta} \bigl
\langle\Pi^\bot e_j,WVu\bigr\rangle,
\]
where
$e_j$ is the $j$th column of $I_{p\times p}$.
Taking $\delta=1/2$, by Proposition \ref{prosphere-covering} we have
$\llvert N_\delta\rrvert\le5^d$.

Now $\Pi^{\perp} \Sigma V = 0$ and so
\[
\bigl\langle\Pi^\perp e_j,W V u \bigr\rangle=
\frac{1}{n} \sum_{i=1}^n \bigl\langle
X_i,\Pi^\perp e_j \bigr\rangle\langle
X_i,V u \rangle
\]
is the sum of independent random variables with mean zero.
By Proposition \ref{proorlicz-cauchy-schwarz}, the summands
satisfy
\begin{eqnarray*}
\bigl\llVert\bigl\langle X_i,\Pi^\perp e_j
\bigr\rangle\langle X_i,V u \rangle\bigr\rrVert_{\psi_1} &
\leq&\bigl\llVert\bigl\langle X_i,\Pi^\perp
e_j \bigr\rangle\bigr\rrVert_{\psi_2} \bigl\llVert\langle
X_i,V u \rangle\bigr\rrVert_{\psi_2}
\\
&=& \bigl\llVert\bigl\langle Z_i,\Sigma^{1/2}
\Pi^\perp e_j \bigr\rangle\bigr\rrVert_{\psi_2}
\bigl\llVert\bigl\langle Z_i,\Sigma^{1/2} V u \bigr\rangle
\bigr\rrVert_{\psi_2}
\\
&\leq&\llVert Z_1 \rrVert_{\psi_2}^2 \bigl
\llVert\Sigma^{1/2} \Pi^{\perp} e_j \bigr\rrVert
_2 \bigl\llVert\Sigma^{1/2} V u \bigr\rrVert_2
\\
&\leq&\llVert Z_1 \rrVert_{\psi_2}^2 \sqrt{
\lambda_1 \lambda_{d+1} }.
\end{eqnarray*}
Recall that $\llVert Z \rrVert_{\psi_2}^2=1$.
Then Bernstein's inequality (Lemma \ref{lembernstein}) implies that
for all $t > 0$
and every $u \in\mathcal{N}_{\delta}$
\begin{eqnarray*}
\Pr\bigl( \bigl\llVert\Pi^\bot W\Pi\bigr\rrVert_{2,\infty}>t
\bigr) &\le& \Pr\Bigl(\max_{1\le j\le p}\max_{u\in N_\delta}
\bigl\langle\Pi^\bot e_j,WVu\bigr\rangle> t/2 \Bigr)
\\
&\le&p 5^d\prob\bigl( \bigl\llvert\bigl\langle\Pi^\perp
e_j,W V u \bigr\rangle\bigr\rrvert> t/2 \bigr)
\\
&\leq& 2 p5^d\exp\biggl( - \frac{t^2/8}{2 \lambda_1\lambda_{d+1}/n + t
\sqrt{\lambda_1\lambda_{d+1}} / n} \biggr).
\end{eqnarray*}
\upqed
\end{pf*}

\subsection{The quadratic terms}
%
%
\begin{lemma}\label{lemlower-quadratic}
Let $\varepsilon\geq0$, $q \in(0,1]$, and
\begin{eqnarray*}
\phi(R_q, \varepsilon) &=& \sup\bigl\{ \bigl\langle\SampleCov- \Sigma,
\Pi^\perp U U^T \Pi^\perp\bigr\rangle\dvtx U \in
\mathbb{V}_{p,d}, \llVert U \rrVert_{2,q}^q \leq
R_q,
\\
&&\hspace*{175.3pt} \bigl\llVert\Pi^\perp U \bigr\rrVert_F \leq\varepsilon
\bigr\}.
\end{eqnarray*}
There exist constants $c > 0$ and $c_1$ such that
for all $x\ge c_1$,
\[
\prob\biggl[\phi(R_q, \varepsilon) \geq cx \llVert Z_1
\rrVert_{\psi_2}^2 \lambda_{d+1} \biggl\{ \varepsilon
\frac{E(R_q,\varepsilon)}{\sqrt{n}} + \frac{E^2(R_q,\varepsilon)}{n}
\biggr\} \biggr]\le2 \exp
\bigl(-x^{2/5} \bigr),
\]
where
\[
E(R_q,\varepsilon) = \E\sup\bigl\{ \langle\mathcal{Z},U \rangle\dvtx U
\in
\Real^{p\times d}, \llVert U \rrVert_{2,q}^q \leq2
R_q, \llVert U \rrVert_F \leq\varepsilon\bigr\}
\]
and $\mathcal{Z}$ is a $p \times d$ matrix with i.i.d. $\normal(0,1)$ entries.
As a consequence, we have, for another constant $c_2$
\[
\E\phi(R_q, \varepsilon) \leq c_2 \llVert Z_1
\rrVert_{\psi_2}^2 \lambda_{d+1} \biggl\{ \varepsilon
\frac{E(R_q,\varepsilon)}{\sqrt{n}} + \frac{E^2(R_q,\varepsilon)}{n}
\biggr\}.
\]

Moreover, we have, for another numerical constant $c'$,
%
%
\begin{equation}
\frac{E(R_q,\varepsilon)}{\sqrt{n}} \leq c'\bigl(R_q^{1/2}
t^{1-q/2} \varepsilon+ R_q t^{2-q}\bigr)
\end{equation}
with $t = \sqrt{\frac{d + \log p}{n}}$.
\end{lemma}
\begin{pf}
The first part follows from Corollary 4.1 of \citet{Vu2012b}. It
remains for us
to prove the ``moreover'' part. By the duality of the $(2,1)$- and
$(2,\infty)$-norms,
\[
\langle\mathcal{Z},U \rangle\leq\llVert\mathcal{Z} \rrVert_{2,\infty}
\llVert U \rrVert_{2,1}
\]
and so
\[
\E(R_q,\varepsilon) \leq\E\llVert\mathcal{Z} \rrVert_{2,\infty}
\sup\bigl\{ \llVert U \rrVert_{2,1}\dvtx U \in\Real^{p\times d},
\llVert U \rrVert_{2,q}^q \leq2 R_q, \llVert U
\rrVert_F \leq\varepsilon\bigr\}.
\]
By (\ref{eqgaussian-row-norm-bound}) and the fact that the Orlicz
$\psi_2$-norm bounds the expectation,
\[
\E\llVert\mathcal{Z} \rrVert_{2,\infty} \leq c' \sqrt{d+
\log p}.
\]
Now $\llVert U \rrVert_{2,1}$ is just the $\ell_1$ norm of
the vector of row-wise
norms of $U$. So we use a standard argument to bound the $\ell_1$ norm
in terms
of the $\ell_2$ and $\ell_q$ norms for $q \in(0,1]$
[e.g., \citet{Raskutti2011}, Lemma 5], and find that for every $t
> 0$
\begin{eqnarray*}
\llVert U \rrVert_{2,1} &\leq& \llVert U \rrVert_{2,q}^{q/2}
\llVert U \rrVert_{2,2} t^{-q/2} + \llVert U \rrVert
_{2,q}^q t^{1-q}
\\
&=& \llVert U \rrVert_{2,q}^{q/2} \llVert U \rrVert
_F t^{-q/2} + \llVert U \rrVert_{2,q}^q
t^{1-q}.
\end{eqnarray*}
Thus,
\[
\sup\bigl\{ \llVert U \rrVert_{2,1}\dvtx U \in\Real^{p\times d},
\llVert U \rrVert_{2,q}^q \leq2 R_q, \llVert U
\rrVert_F \leq\varepsilon\bigr\} \leq R_q^{1/2}
t^{-q/2} + R_q t^{1-q}.
\]
Letting $t = \E\llVert\mathcal{Z} \rrVert_{2,\infty} /
\sqrt{n}$,
and combining the above inequalities completes the proof.
\end{pf}


%
%
\begin{lemma}\label{lemupper-quadratic}
There exists a constant $c > 0$ such that
\[
\bigl\llVert{\bigl\llVert\Pi(\SampleCov-\Sigma)\Pi\bigr\rrVert_2}
\bigr\rrVert_{\psi_1} \leq c \llVert Z_1 \rrVert
_{\psi_2}^2 \lambda_1 ( \sqrt{d/n} + d/n ).
\]
\end{lemma}
\begin{pf}
Let $\mathcal{N}_\delta$ be a minimal $\delta$-net of $\Sphere
_2^{d-1}$ for
some $\delta\in(0,1)$ to be chosen later. Then
\[
\bigl\llVert\Pi(\SampleCov-\Sigma)\Pi\bigr\rrVert_2 = \bigl\llVert
V^T (\SampleCov-\Sigma)V \bigr\rrVert_2 \leq(1-2
\delta)^{-1} \max_{u \in\mathcal{N}_{\delta}} \bigl\llvert\bigl\langle V
u,(\SampleCov-\Sigma) V u \bigr\rangle\bigr\rrvert.
\]
Using a similar argument as in the proof of Lemma \ref{lemcrosstermtail},
for all $t > 0$ and every $u \in\mathcal{N}_\delta$
\[
\prob\bigl( \bigl\llvert\bigl\langle V u,(\SampleCov- \Sigma) V u
\bigr\rangle
\bigr\rrvert> t \bigr) \leq2 \exp\biggl( - \frac{t^2/2}{2 \sigma^2/n
+ t \sigma/ n} \biggr),
\]
where $\sigma= 2 \llVert Z_1 \rrVert_{\psi_2}^2 \lambda_1$.
Then Lemma \ref{lembernstein-maximal} implies that
\begin{eqnarray*}
\bigl\Vert{ \bigl\llVert\Pi(\SampleCov-\Sigma)\Pi\bigr\Vert_2\bigr\rrVert
}_{\psi_1} &\leq& (1 - 2 \delta)^{-1} \Bigl\Vert{ \max
_{u \in\mathcal{N}_{\delta}} \bigl\llvert\bigl\langle V u,(\SampleCov
-\Sigma) V u
\bigr\rangle\bigr\rrvert}\Bigr\Vert_{\psi_1}
\\
&\leq& (1 - 2 \delta)^{-1} C \sigma\biggl( \sqrt{\frac{ \log(1 +
\llvert\mathcal{N}_{\delta}\rrvert
) }{ n }} +
\frac{ \log(1 + \llvert\mathcal{N}_{\delta}\rrvert) }{ n } \biggr),
\end{eqnarray*}
where $C > 0$ is a constant.
Choosing $\delta= 1/3$ and applying Proposition \ref
{prosphere-covering} yields
$\llvert\mathcal{N}_{\delta}\rrvert\leq7^d$ and
\[
\log\bigl(1 + \llvert\mathcal{N}_{\delta}\rrvert\bigr) \leq\log(8)
\log(d).
\]
Thus,
\[
\bigl\llVert{ \bigl\llVert\Pi(\SampleCov-\Sigma)\Pi\bigr\rrVert_2
}\bigr\rrVert_{\psi_1} \leq7 C \sigma( \sqrt{d/n} + d/n ).
\]
\upqed
\end{pf}
\end{appendix}

\section*{Acknowledgments}

The research reported in this article was completed while V. Q. Vu was visiting
the Department of Statistics at Carnegie Mellon University. He thanks
them for
their hospitality and support. We also thank the referees for their
helpful comments.


%

\printaddresses


\begin{thebibliography}{38}

\bibitem[\protect\citeauthoryear{Amini and Wainwright}{2009}]{Amini2009}
%
\begin{barticle}[mr]
\bauthor{\bsnm{Amini},~\bfnm{Arash~A.}\binits{A.~A.}} \AND
\bauthor{\bsnm{Wainwright},~\bfnm{Martin~J.}\binits{M.~J.}}
(\byear{2009}).
\btitle{High-dimensional analysis of semidefinite relaxations for sparse
principal components}.
\bjournal{Ann. Statist.}
\bvolume{37}
\bpages{2877--2921}.
\bid{doi={10.1214/08-AOS664}, issn={0090-5364}, mr={2541450}}
\bptok{imsref}%
\end{barticle}
%
\endbibitem

\bibitem[\protect\citeauthoryear{Bhatia}{1997}]{Bhatia}
%
\begin{bbook}[mr]
\bauthor{\bsnm{Bhatia},~\bfnm{Rajendra}\binits{R.}}
(\byear{1997}).
\btitle{Matrix Analysis}.
\bseries{Graduate Texts in Mathematics}
\bvolume{169}.
\bpublisher{Springer}, \blocation{New York}.
\bid{doi={10.1007/978-1-4612-0653-8}, mr={1477662}}
\bptok{imsref}%
\end{bbook}
%
\endbibitem

\bibitem[\protect\citeauthoryear{Birnbaum et~al.}{2013}]{Birnbaum2012}
\begin{barticle}[mr]
\bauthor{\bsnm{Birnbaum},~\bfnm{Aharon}\binits{A.}},
  \bauthor{\bsnm{Johnstone},~\bfnm{Iain~M.}\binits{I.~M.}},
  \bauthor{\bsnm{Nadler},~\bfnm{Boaz}\binits{B.}} \AND
  \bauthor{\bsnm{Paul},~\bfnm{Debashis}\binits{D.}}
(\byear{2013}).
\btitle{Minimax bounds for sparse {PCA} with noisy high-dimensional data}.
\bjournal{Ann. Statist.}
\bvolume{41}
\bpages{1055--1084}.
\bid{issn={0090-5364}, mr={3113803}}
\bptok{imsref}%
\end{barticle}
\endbibitem


\bibitem[\protect\citeauthoryear{Chen, Zou and Cook}{2010}]{Chen2010}
%
\begin{barticle}[mr]
\bauthor{\bsnm{Chen},~\bfnm{Xin}\binits{X.}},
\bauthor{\bsnm{Zou},~\bfnm{Changliang}\binits{C.}} \AND
\bauthor{\bsnm{Cook},~\bfnm{R.~Dennis}\binits{R.~D.}}
(\byear{2010}).
\btitle{Coordinate-independent sparse sufficient dimension reduction and
variable selection}.
\bjournal{Ann. Statist.}
\bvolume{38}
\bpages{3696--3723}.
\bid{doi={10.1214/10-AOS826}, issn={0090-5364}, mr={2766865}}
\bptok{imsref}%
\end{barticle}
%
\endbibitem

\bibitem[\protect\citeauthoryear{Chikuse}{2003}]{Chikuse2003}
%
\begin{bbook}[mr]
\bauthor{\bsnm{Chikuse},~\bfnm{Yasuko}\binits{Y.}}
(\byear{2003}).
\btitle{Statistics on Special Manifolds}.
\bseries{Lecture Notes in Statistics}
\bvolume{174}.
\bpublisher{Springer}, \blocation{New York}.
\bid{mr={1960435}}
\bptok{imsref}%
\end{bbook}
%
\endbibitem

\bibitem[\protect\citeauthoryear{d'Aspremont et~al.}{2007}]{dAspremont2007}
%
\begin{barticle}[mr]
\bauthor{\bsnm{d'Aspremont},~\bfnm{Alexandre}\binits{A.}},
\bauthor{\bsnm{El~Ghaoui},~\bfnm{Laurent}\binits{L.}},
\bauthor{\bsnm{Jordan},~\bfnm{Michael~I.}\binits{M.~I.}} \AND
\bauthor{\bsnm{Lanckriet},~\bfnm{Gert R.~G.}\binits{G.~R.~G.}}
(\byear{2007}).
\btitle{A direct formulation for sparse {PCA} using semidefinite programming}.
\bjournal{SIAM Rev.}
\bvolume{49}
\bpages{434--448 (electronic)}.
\bid{doi={10.1137/050645506}, issn={0036-1445}, mr={2353806}}
\bptok{imsref}%
\end{barticle}
%
\endbibitem

\bibitem[\protect\citeauthoryear{Donoho and Johnstone}{1994}]{Donoho1994}
%
\begin{barticle}[author]
\bauthor{\bsnm{Donoho},~\bfnm{David~L}\binits{D.~L.}} \AND
\bauthor{\bsnm{Johnstone},~\bfnm{Iain~M}\binits{I.~M.}}
(\byear{1994}).
\btitle{Minimax risk over $l_p$-balls for $l_q$-error}.
\bjournal{Probab. Theory Related Fields}
\bvolume{99}
\bpages{277--303}.
\bptok{imsref}%
\end{barticle}
%
\endbibitem

\bibitem[\protect\citeauthoryear{Edelman, Arias and Smith}{1999}]{Edelman1998}
%
\begin{barticle}[mr]
\bauthor{\bsnm{Edelman},~\bfnm{Alan}\binits{A.}},
\bauthor{\bsnm{Arias},~\bfnm{Tom{\'a}s~A.}\binits{T.~A.}} \AND
\bauthor{\bsnm{Smith},~\bfnm{Steven~T.}\binits{S.~T.}}
(\byear{1999}).
\btitle{The geometry of algorithms with orthogonality constraints}.
\bjournal{SIAM J. Matrix Anal. Appl.}
\bvolume{20}
\bpages{303--353}.
\bid{doi={10.1137/S0895479895290954}, issn={0895-4798}, mr={1646856}}
\bptnote{check year}%
\bptok{imsref}%
\end{barticle}
%
\endbibitem

\bibitem[\protect\citeauthoryear{Gilbert}{1952}]{Gilbert52}
%
\begin{barticle}[author]
\bauthor{\bsnm{Gilbert},~\bfnm{E.~N.}\binits{E.~N.}}
(\byear{1952}).
\btitle{A comparison of signalling alphabets}.
\bjournal{Bell System Technical Journal}
\bvolume{31}
\bpages{504--522}.
\bptok{imsref}%
\end{barticle}
%
\endbibitem

\bibitem[\protect\citeauthoryear{Hotelling}{1933}]{Hotelling1933}
%
\begin{barticle}[author]
\bauthor{\bsnm{Hotelling},~\bfnm{Harold}\binits{H.}}
(\byear{1933}).
\btitle{Analysis of a complex of statistical variables into principal
components}.
\bjournal{Journal of Educational Psychology}
\bvolume{24}
\bpages{498--520}.
\bptok{imsref}%
\end{barticle}
%
\endbibitem

\bibitem[\protect\citeauthoryear{Johnstone and Lu}{2009}]{Johnstone2009}
%
\begin{barticle}[mr]
\bauthor{\bsnm{Johnstone},~\bfnm{Iain~M.}\binits{I.~M.}} \AND
\bauthor{\bsnm{Lu},~\bfnm{Arthur~Yu}\binits{A.~Y.}}
(\byear{2009}).
\btitle{On consistency and sparsity for principal components analysis
in high
dimensions}.
\bjournal{J. Amer. Statist. Assoc.}
\bvolume{104}
\bpages{682--693}.
\bid{doi={10.1198/jasa.2009.0121}, issn={0162-1459}, mr={2751448}}
\bptok{imsref}%
\end{barticle}
%
\endbibitem

\bibitem[\protect\citeauthoryear{Jolliffe, Trendafilov and
Uddin}{2003}]{Jolliffe2003}
%
\begin{barticle}[mr]
\bauthor{\bsnm{Jolliffe},~\bfnm{Ian~T.}\binits{I.~T.}},
\bauthor{\bsnm{Trendafilov},~\bfnm{Nickolay~T.}\binits{N.~T.}} \AND
\bauthor{\bsnm{Uddin},~\bfnm{Mudassir}\binits{M.}}
(\byear{2003}).
\btitle{A modified principal component technique based on the {LASSO}}.
\bjournal{J. Comput. Graph. Statist.}
\bvolume{12}
\bpages{531--547}.
\bid{doi={10.1198/1061860032148}, issn={1061-8600}, mr={2002634}}
\bptok{imsref}%
\end{barticle}
%
\endbibitem

\bibitem[\protect\citeauthoryear{Journ{\'e}e et~al.}{2010}]{Journee2010}
%
\begin{barticle}[mr]
\bauthor{\bsnm{Journ{\'e}e},~\bfnm{Michel}\binits{M.}},
\bauthor{\bsnm{Nesterov},~\bfnm{Yurii}\binits{Y.}},
\bauthor{\bsnm{Richt{\'a}rik},~\bfnm{Peter}\binits{P.}} \AND
\bauthor{\bsnm{Sepulchre},~\bfnm{Rodolphe}\binits{R.}}
(\byear{2010}).
\btitle{Generalized power method for sparse principal component analysis}.
\bjournal{J. Mach. Learn. Res.}
\bvolume{11}
\bpages{517--553}.
\bid{issn={1532-4435}, mr={2600619}}
\bptok{imsref}%
\end{barticle}
%
\endbibitem

\bibitem[\protect\citeauthoryear{Ledoux}{2001}]{Ledoux}
%
\begin{bbook}[mr]
\bauthor{\bsnm{Ledoux},~\bfnm{Michel}\binits{M.}}
(\byear{2001}).
\btitle{The Concentration of Measure Phenomenon}.
\bseries{Mathematical Surveys and Monographs}
\bvolume{89}.
\bpublisher{Amer. Math. Soc.}, \blocation{Providence, RI}.
\bid{mr={1849347}}
\bptok{imsref}%
\end{bbook}
%
\endbibitem

\bibitem[\protect\citeauthoryear{Lounici}{2013}]{Lounici2012}
\begin{bincollection}[author]
\bauthor{\bsnm{Lounici},~\bfnm{Karim}\binits{K.}}
(\byear{2013}).
\btitle{Sparse principal component analysis with missing observations}.
In  \bbooktitle{High Dimensional Probability VI}
(\beditor{\bfnm{C.}\binits{C.}~\bsnm{Houdr\'e}},
  \beditor{\bfnm{D.~M.}\binits{D.~M.}~\bsnm{Mason}},
  \beditor{\bfnm{J.}\binits{J.}~\bsnm{Rosi\'nski}} \AND
  \beditor{\bfnm{J.~A.}\binits{J.~A.}~\bsnm{Wellner}}, eds.).
\bseries{Progress in Probability}
\bvolume{66}
\bpages{327--356}.
\bpublisher{Springer}, \blocation{Basel}.
\bptok{imsref}%
\end{bincollection}
%
\endbibitem


\bibitem[\protect\citeauthoryear{Ma}{2013}]{Ma2011}
\begin{barticle}[mr]
\bauthor{\bsnm{Ma},~\bfnm{Zongming}\binits{Z.}}
(\byear{2013}).
\btitle{Sparse principal component analysis and iterative thresholding}.
\bjournal{Ann. Statist.}
\bvolume{41}
\bpages{772--801}.
\bid{doi={10.1214/13-AOS1097}, issn={0090-5364}, mr={3099121}}
\bptok{imsref}%
\end{barticle}
\endbibitem


\bibitem[\protect\citeauthoryear{Massart}{2007}]{Massart2007}
%
\begin{bbook}[mr]
\bauthor{\bsnm{Massart},~\bfnm{Pascal}\binits{P.}}
(\byear{2007}).
\btitle{Concentration Inequalities and Model Selection}.
\bseries{Lecture Notes in Math.}
\bvolume{1896}.
\bpublisher{Springer}, \blocation{Berlin}.
\bid{mr={2319879}}
\bptok{imsref}%
\end{bbook}
%
\endbibitem

\bibitem[\protect\citeauthoryear{Mendelson}{2010}]{Mendelson2010}
%
\begin{barticle}[mr]
\bauthor{\bsnm{Mendelson},~\bfnm{Shahar}\binits{S.}}
(\byear{2010}).
\btitle{Empirical processes with a bounded {$\psi\sb1$} diameter}.
\bjournal{Geom. Funct. Anal.}
\bvolume{20}
\bpages{988--1027}.
\bid{doi={10.1007/s00039-010-0084-5}, issn={1016-443X}, mr={2729283}}
\bptok{imsref}%
\end{barticle}
%
\endbibitem

\bibitem[\protect\citeauthoryear{Nadler}{2008}]{Nadler2008}
%
\begin{barticle}[mr]
\bauthor{\bsnm{Nadler},~\bfnm{Boaz}\binits{B.}}
(\byear{2008}).
\btitle{Finite sample approximation results for principal component
analysis: A
matrix perturbation approach}.
\bjournal{Ann. Statist.}
\bvolume{36}
\bpages{2791--2817}.
\bid{doi={10.1214/08-AOS618}, issn={0090-5364}, mr={2485013}}
\bptok{imsref}%
\end{barticle}
%
\endbibitem

\bibitem[\protect\citeauthoryear{Negahban et~al.}{2012}]{Negahban2012}
%
\begin{barticle}[mr]
\bauthor{\bsnm{Negahban},~\bfnm{Sahand~N.}\binits{S.~N.}},
\bauthor{\bsnm{Ravikumar},~\bfnm{Pradeep}\binits{P.}},
\bauthor{\bsnm{Wainwright},~\bfnm{Martin~J.}\binits{M.~J.}} \AND
\bauthor{\bsnm{Yu},~\bfnm{Bin}\binits{B.}}
(\byear{2012}).
\btitle{A unified framework for high-dimensional analysis of {$M$}-estimators
with decomposable regularizers}.
\bjournal{Statist. Sci.}
\bvolume{27}
\bpages{538--557}.
\bid{doi={10.1214/12-STS400}, issn={0883-4237}, mr={3025133}}
\bptok{imsref}%
\end{barticle}
%
\endbibitem

\bibitem[\protect\citeauthoryear{Pajor}{1998}]{Pajor1998}
%
\begin{bincollection}[mr]
\bauthor{\bsnm{Pajor},~\bfnm{Alain}\binits{A.}}
(\byear{1998}).
\btitle{Metric entropy of the {G}rassmann manifold}.
In \bbooktitle{Convex Geometric Analysis ({B}erkeley, {CA}, 1996)}.
\bseries{Mathematical Sciences Research Institute Publications}
\bvolume{34}
\bpages{181--188}.
\bpublisher{Cambridge Univ. Press}, \blocation{Cambridge}.
\bid{mr={1665590}}
\bptnote{check year}%
\bptok{imsref}%
\end{bincollection}
%
\endbibitem

\bibitem[\protect\citeauthoryear{Paul}{2007}]{Paul2007}
%
\begin{barticle}[mr]
\bauthor{\bsnm{Paul},~\bfnm{Debashis}\binits{D.}}
(\byear{2007}).
\btitle{Asymptotics of sample eigenstructure for a large dimensional spiked
covariance model}.
\bjournal{Statist. Sinica}
\bvolume{17}
\bpages{1617--1642}.
\bid{issn={1017-0405}, mr={2399865}}
\bptok{imsref}%
\end{barticle}
%
\endbibitem

\bibitem[\protect\citeauthoryear{Pearson}{1901}]{Pearson1901}
%
\begin{barticle}[author]
\bauthor{\bsnm{Pearson},~\bfnm{Karl}\binits{K.}}
(\byear{1901}).
\btitle{On lines and planes of closest fit to systems of points in space}.
\bjournal{Philosophical Magazine}
\bvolume{2}
\bpages{559--572}.
\bptok{imsref}%
\end{barticle}
%
\endbibitem

\bibitem[\protect\citeauthoryear{Raskutti, Wainwright and
  Yu}{2011}]{Raskutti2011}
\begin{barticle}[mr]
\bauthor{\bsnm{Raskutti},~\bfnm{Garvesh}\binits{G.}},
  \bauthor{\bsnm{Wainwright},~\bfnm{Martin~J.}\binits{M.~J.}} \AND
  \bauthor{\bsnm{Yu},~\bfnm{Bin}\binits{B.}}
(\byear{2011}).
\btitle{Minimax rates of estimation for high-dimensional linear regression over
  {$\ell\sb q$}-balls}.
\bjournal{IEEE Trans. Inform. Theory}
\bvolume{57}
\bpages{6976--6994}.
\bid{doi={10.1109/TIT.2011.2165799}, issn={0018-9448}, mr={2882274}}
\bptok{imsref}%
\end{barticle}
\endbibitem


\bibitem[\protect\citeauthoryear{Shen and Huang}{2008}]{Shen2008}
%
\begin{barticle}[mr]
\bauthor{\bsnm{Shen},~\bfnm{Haipeng}\binits{H.}} \AND
\bauthor{\bsnm{Huang},~\bfnm{Jianhua~Z.}\binits{J.~Z.}}
(\byear{2008}).
\btitle{Sparse principal component analysis via regularized low rank matrix
approximation}.
\bjournal{J. Multivariate Anal.}
\bvolume{99}
\bpages{1015--1034}.
\bid{doi={10.1016/j.jmva.2007.06.007}, issn={0047-259X}, mr={2419336}}
\bptok{imsref}%
\end{barticle}
%
\endbibitem

\bibitem[\protect\citeauthoryear{Shen, Shen and Marron}{2013}]{Shen2011}
\begin{barticle}[mr]
\bauthor{\bsnm{Shen},~\bfnm{Dan}\binits{D.}},
  \bauthor{\bsnm{Shen},~\bfnm{Haipeng}\binits{H.}} \AND
  \bauthor{\bsnm{Marron},~\bfnm{J.~S.}\binits{J.~S.}}
(\byear{2013}).
\btitle{Consistency of sparse {PCA} in high dimension, low sample size
  contexts}.
\bjournal{J. Multivariate Anal.}
\bvolume{115}
\bpages{317--333}.
\bid{doi={10.1016/j.jmva.2012.10.007}, issn={0047-259X}, mr={3004561}}
\bptok{imsref}%
\end{barticle}
\endbibitem


\bibitem[\protect\citeauthoryear{Stewart and Sun}{1990}]{StewartAndSun}
%
\begin{bbook}[mr]
\bauthor{\bsnm{Stewart},~\bfnm{G.~W.}\binits{G.~W.}} \AND
\bauthor{\bsnm{Sun},~\bfnm{Ji~Guang}\binits{J.~G.}}
(\byear{1990}).
\btitle{Matrix Perturbation Theory}.
\bpublisher{Academic Press}, \blocation{Boston, MA}.
\bid{mr={1061154}}
\bptok{imsref}%
\end{bbook}
%
\endbibitem

\bibitem[\protect\citeauthoryear{Tibshirani}{1996}]{Tibshirani1996}
%
\begin{barticle}[mr]
\bauthor{\bsnm{Tibshirani},~\bfnm{Robert}\binits{R.}}
(\byear{1996}).
\btitle{Regression shrinkage and selection via the lasso}.
\bjournal{J. R. Stat. Soc. Ser. B Stat. Methodol.}
\bvolume{58}
\bpages{267--288}.
\bid{issn={0035-9246}, mr={1379242}}
\bptok{imsref}%
\end{barticle}
%
\endbibitem

\bibitem[\protect\citeauthoryear{van~der Vaart and
Wellner}{1996}]{vanderVaartAndWellner}
%
\begin{bbook}[mr]
\bauthor{\bparticle{van~der} \bsnm{Vaart},~\bfnm{Aad~W.}\binits{A.~W.}}
\AND
\bauthor{\bsnm{Wellner},~\bfnm{Jon~A.}\binits{J.~A.}}
(\byear{1996}).
\btitle{Weak Convergence and Empirical Processes}.
\bpublisher{Springer}, \blocation{New York}.
\bid{mr={1385671}}
\bptok{imsref}%
\end{bbook}
%
\endbibitem

\bibitem[\protect\citeauthoryear{Var{\v{s}}amov}{1957}]{Varshamov57}
%
\begin{barticle}[mr]
\bauthor{\bsnm{Var{\v{s}}amov},~\bfnm{R.~R.}\binits{R.~R.}}
(\byear{1957}).
\btitle{The evaluation of signals in codes with correction of errors}.
\bjournal{Dokl. Akad. Nauk SSSR (N.S)}
\bvolume{117}
\bpages{739--741}.
\bid{issn={0002-3264}, mr={0095090}}
\bptnote{check related}%
\bptok{imsref}%
\end{barticle}
%
\endbibitem

\bibitem[\protect\citeauthoryear{Vu and Lei}{2012a}]{Vu2012}
%
\begin{binproceedings}[author]
\bauthor{\bsnm{Vu},~\bfnm{Vincent~Q.}\binits{V.~Q.}} \AND
\bauthor{\bsnm{Lei},~\bfnm{Jing}\binits{J.}}
(\byear{2012}a).
\btitle{Minimax rates of estimation for sparse PCA in high dimensions}.
In \bbooktitle{Proceedings of the 15th International Conference on Artificial
Intelligence and Statistics \mbox{(AISTATS)}.
JMLR Workshop and Conference Proceedings Volume 22.}
\bptok{imsref}%
\end{binproceedings}
%
\endbibitem

\bibitem[\protect\citeauthoryear{Vu and Lei}{2012b}]{Vu2012b}
%
\begin{bmisc}[author]
\bauthor{\bsnm{Vu},~\bfnm{Vincent~Q.}\binits{V.~Q.}} \AND
\bauthor{\bsnm{Lei},~\bfnm{Jing}\binits{J.}}
(\byear{2012}b).
\bhowpublished{Squared-norm empirical process in {B}anach space.
Available at arXiv:\arxivurl{1312.1005}.}
\bptok{imsref}%
\end{bmisc}
%
\endbibitem

\bibitem[\protect\citeauthoryear{Witten, Tibshirani and
Hastie}{2009}]{Witten2009}
%
\begin{barticle}[pbm]
\bauthor{\bsnm{Witten},~\bfnm{Daniela~M.}\binits{D.~M.}},
\bauthor{\bsnm{Tibshirani},~\bfnm{Robert}\binits{R.}} \AND
\bauthor{\bsnm{Hastie},~\bfnm{Trevor}\binits{T.}}
(\byear{2009}).
\btitle{A penalized matrix decomposition, with applications to sparse principal
components and canonical correlation analysis}.
\bjournal{Biostatistics}
\bvolume{10}
\bpages{515--534}.
\bid{doi={10.1093/biostatistics/kxp008}, issn={1468-4357}, pii={kxp008},
pmcid={2697346}, pmid={19377034}}
\bptok{imsref}%
\end{barticle}
%
\endbibitem

\bibitem[\protect\citeauthoryear{Yu}{1997}]{Yu1997}
%
\begin{bincollection}[author]
\bauthor{\bsnm{Yu},~\bfnm{Bin}\binits{B.}}
(\byear{1997}).
\btitle{{Assouad}, {Fano}, and {Le Cam}}.
In \bbooktitle{Festschrift for {Lucien Le Cam}}
(\beditor{\bfnm{D.}\binits{D.}~\bsnm{Pollard}},
  \beditor{\bfnm{E.}\binits{E.}~\bsnm{Torgersen}} \AND
  \beditor{\bfnm{G.~L.}\binits{G.~L.}~\bsnm{Yang}}, eds.)
\bpages{423--435}.
\bpublisher{Springer}, \blocation{New York}.
\bptok{imsref}%
\end{bincollection}
%
\endbibitem

\bibitem[\protect\citeauthoryear{Yuan and Lin}{2006}]{Yuan2006}
\begin{barticle}[mr]
\bauthor{\bsnm{Yuan},~\bfnm{Ming}\binits{M.}} \AND
  \bauthor{\bsnm{Lin},~\bfnm{Yi}\binits{Y.}}
(\byear{2006}).
\btitle{Model selection and estimation in regression with grouped variables}.
\bjournal{J. R. Stat. Soc. Ser. B Stat. Methodol.}
\bvolume{68}
\bpages{49--67}.
\bid{doi={10.1111/j.1467-9868.2005.00532.x}, issn={1369-7412}, mr={2212574}}
\bptok{imsref}%
\end{barticle}
\endbibitem


\bibitem[\protect\citeauthoryear{Zhao, Rocha and Yu}{2009}]{Zhao2009}
%
\begin{barticle}[mr]
\bauthor{\bsnm{Zhao},~\bfnm{Peng}\binits{P.}},
\bauthor{\bsnm{Rocha},~\bfnm{Guilherme}\binits{G.}} \AND
\bauthor{\bsnm{Yu},~\bfnm{Bin}\binits{B.}}
(\byear{2009}).
\btitle{The composite absolute penalties family for grouped and hierarchical
variable selection}.
\bjournal{Ann. Statist.}
\bvolume{37}
\bpages{3468--3497}.
\bid{doi={10.1214/07-AOS584}, issn={0090-5364}, mr={2549566}}
\bptok{imsref}%
\end{barticle}
%
\endbibitem

\bibitem[\protect\citeauthoryear{Zou, Hastie and Tibshirani}{2006}]{Zou2006}
%
\begin{barticle}[mr]
\bauthor{\bsnm{Zou},~\bfnm{Hui}\binits{H.}},
\bauthor{\bsnm{Hastie},~\bfnm{Trevor}\binits{T.}} \AND
\bauthor{\bsnm{Tibshirani},~\bfnm{Robert}\binits{R.}}
(\byear{2006}).
\btitle{Sparse principal component analysis}.
\bjournal{J.~Comput. Graph. Statist.}
\bvolume{15}
\bpages{265--286}.
\bid{doi={10.1198/106186006X113430}, issn={1061-8600}, mr={2252527}}
\bptok{imsref}%
\end{barticle}
%
\endbibitem

\end{thebibliography}
\end{document}